\input amstex
\documentstyle{amsppt}
\magnification1200

\def\A{{\Cal A}}
\def\al{\alpha}

\def\B{{\Cal B}}
\def\ba{{\bar a}}
\def\bb{{\bar b}}
\def\bi{{\bar\imath}}
\def\bj{{\bar\jmath}}

\def\bx{{\boxed{\phantom{\square}}\kern-.4pt}}

\def\C{{\Cal C}}
\def\CC{{\Bbb C}}
\def\cj{{\tilde\jmath}}
\def\Chi{{\operatorname{X}}}
\def\CM{\CC^{M|M}}
\def\CMN{\CC^{N|N}\!\ot\CC^{M|M}}
\def\CN{\CC^{N|N}}
\def\crc{{\raise.24ex\hbox{$\sssize\kern.1em\circ\kern.1em$}}}

\def\d{\partial}
\def\de{\delta}

\def\dim{{\operatorname{dim\,}}}
\def\DMN{{\Cal{D}}}

\def\End{\operatorname{End}\hskip1pt}
\def\EndCN{\End(\Bbb C^{N|N})}
\def\enddemos{{\ $\square$\enddemo}}
\def\ep{\varepsilon}

\def\ga{\gamma}
\def\Ga{\Gamma}
\def\Gaz{{\Ga(t_1\lc t_m)}}
\def\ge{\geqslant}

\def\glN{{\frak{gl}_{N|N}}}

\def\hs{{\hskip1pt}}

\def\I{{\Cal I}}
\def\id{{\operatorname{id}}}

\def\J{{\,\Cal J}}

\def\la{\lambda}

\def\La{\Lambda}

\def\le{\leqslant}

\def\om{\omega}
\def\ot{\otimes}

\def\PDMN{{\Cal{PD}}}
\def\PMN{{\Cal{P}}}
\def\phi{\varphi}
\def\phik{{\phi_{ijk}(u,v,w)}}
\def\Psiz{{\Psi_{\la}(t_1\lc t_n)}}

\def\qM{{\frak{q}_M}}
\def\qN{{\frak{q}_N}}

\def\Sen{H_n}

\def\The{\Theta}
\def\Thez{{\The_\la(t_1\lc t_n)}}

\def\tr{{\operatorname{str}}}
\def\ts{{\hskip1pt}}

\def\UM{{\operatorname{U}(\qM)}}
\def\UN{{\operatorname{U}(\qN)}}

\def\Up{\Upsilon}
\def\Upz{{\Up_\la(t_1\lc t_n)}}

\def\Xiu{{\operatorname{X}(u)}}

\def\Z{{\Cal T}}
\def\ZM{{\operatorname{Z}(\qM)}}
\def\ZN{{\operatorname{Z}(\qN)}}
\def\ZZ{{\Bbb Z}}

\expandafter\ifx\csname bethe.def\endcsname\relax \else\endinput\fi
\expandafter\edef\csname bethe.def\endcsname{%
 \catcode`\noexpand\@=\the\catcode`\@\space}
\catcode`\@=11

\mathsurround 1.6pt

\def\hcor#1{\advance\hoffset by #1}
\def\vcor#1{\advance\voffset by #1}
\let\bls\baselineskip  \let\ignore\ignorespaces
\def\vsk#1>{\vskip#1\bls} \let\adv\advance 
\def\vv#1>{\vadjust{\vsk#1>}\ignore} \def\vvv#1>{\vadjust{\vskip#1}\ignore}
\def\vvn#1>{\vadjust{\nobreak\vsk#1>\nobreak}\ignore}
\def\vvvn#1>{\vadjust{\nobreak\vskip#1\nobreak}\ignore}
\def\setnormalbls{\edef\normalbls{\bls\the\bls}}
\def\setmaths{\edef\maths{\mathsurround\the\mathsurround}}

 \let\nt\noindent 
\def\nn#1>{\noalign{\vskip #1pt}} \def\NN#1>{\openup#1pt}
 
\let\Sum\sum \def\sum{\Sum\limits} 
\let\Prod\prod \def\prod{\Prod\limits} \let\Int\int \def\int{\Int\limits}

\let\=\m@th \def\&{.\kern.1em} \def\>{\!\;} \def\:{\!\!\;}

\ifx\plainfootnote\undefined \let\plainfootnote\footnote \fi
\expandafter\ifx\csname amsppt.sty\endcsname\relax
 
\else \fi

\newbox\s@ctb@x
\def\s@ct#1 #2\par{\removelastskip\vsk>
 \vtop{\bf\setbox\s@ctb@x\hbox{#1} \parindent\wd\s@ctb@x
 \ifdim\parindent>0pt\adv\parindent.5em\fi\item{#1}#2\strut}%
 \nointerlineskip\nobreak\vtop{\strut}\nobreak\vsk-.4>\nobreak}

\newbox\t@stb@x
\def\gadv{\global\advance} \def\gad#1{\gadv#1 1} 
\def\l@b@l#1#2{\def\n@@{\csname #2no\endcsname}%
 \if *#1\gad\n@@ \expandafter\xdef\csname @#1@#2@\endcsname{\the\Sno.\the\n@@}%
 \else\expandafter\ifx\csname @#1@#2@\endcsname\relax\gad\n@@
 \expandafter\xdef\csname @#1@#2@\endcsname{\the\Sno.\the\n@@}\fi\fi}
\def\l@bel#1#2{\l@b@l{#1}{#2}\?#1@#2?}
\def\?#1?{\csname @#1@\endcsname}
\def\[#1]{\def\n@xt@{\ifx\t@st *\def\n@xt####1{{\setbox\t@stb@x\hbox{\?#1@F?}%
 \ifnum\wd\t@stb@x=0 {\bf???}\else\?#1@F?\fi}}\else
 \def\n@xt{{\setbox\t@stb@x\hbox{\?#1@L?}\ifnum\wd\t@stb@x=0 {\bf???}\else
 \?#1@L?\fi}}\fi\n@xt}\futurelet\t@st\n@xt@}
\def\(#1){{\rm\setbox\t@stb@x\hbox{\?#1@F?}\ifnum\wd\t@stb@x=0 ({\bf???})\else
 (\?#1@F?)\fi}}
\def\dff{\expandafter\d@f} \def\d@f{\expandafter\def}
\def\edff{\expandafter\ed@f} \def\ed@f{\expandafter\edef}

\newcount\Sno \newcount\Lno \newcount\Fno
\def\Section#1{\gadno\Fno=0\Lno=0\s@ct{\the\Sno.} {#1}\par} \let\Sect\Section
\def\section#1{\gad\Sno\Fno=0\Lno=0\s@ct{} {#1}\par} \let\sect\section
\def\l@F#1{\l@bel{#1}F} \def\<#1>{\l@b@l{#1}F} \def\l@L#1{\l@bel{#1}L}
\def\Tag#1{\tag\l@F{#1}} \def\Tagg#1{\tag"\llap{\rm(\l@F{#1})}"}
\def\Th#1{Theorem \l@L{#1}} \def\Lm#1{Lemma \l@L{#1}}
\def\Prop#1{Proposition \l@L{#1}}
\def\Cr#1{Corollary \l@L{#1}} \def\Cj#1{Conjecture \l@L{#1}}
 
\def\Proof#1.{\demo{\it Proof #1}}

 \def\setparindent{\edef\Parindent{\the\parindent}}
\def\Appendix{\Sno=64\let\p@r@\z@ 
 \def\Section##1{\gad\Sno\Fno=0\Lno=0 \s@ct{} \hskip\p@r@ Appendix
\\the\Sno
  \if *##1\relax\else {.\enspace##1}\fi\par} \let\Sect\Section
 \def\section##1{\gad\Sno\Fno=0\Lno=0 \s@ct{} \hskip\p@r@ Appendix%
  \if *##1\relax\else {.\enspace##1}\fi\par} \let\sect\section
 \def\l@b@l##1##2{\def\n@@{\csname ##2no\endcsname}%
 \if *##1\gad\n@@
 \expandafter\xdef\csname @##1@##2@\endcsname{\char\the\Sno.\the\n@@}%
 \else\expandafter\ifx\csname @##1@##2@\endcsname\relax\gad\n@@
 \expandafter\xdef\csname @##1@##2@\endcsname{\char\the\Sno.\the\n@@}\fi\fi}}

 \let\x\times 
 
\let\le\leqslant \let\ge\geqslant
  \let\8\infty

\let\=\m@th  \def\_#1{_{\rlap{$\ssize#1$}}}

\def\lc{{,\ldots\hskip-0.1pt,\hskip1pt}}

\def\E(#1){\mathop{\hbox{\rm End}\,}(#1)} 

\def\1{^{-1}} \def\vst#1{{\lower2.1pt\hbox{$\bigr|_{#1}$}}}

\let\logo@\relax
\let\m@k@h@@d\makeheadline \let\m@k@f@@t\makefootline
\def\makeheadline{\ifnum\pageno=1\headline={\hfil}\fi\m@k@h@@d}
\def\makefootline{\ifnum\pageno=1\footline={\hfil}\fi\m@k@f@@t}


\font\bigbf=cmbx10 scaled 1200
\centerline{\bigbf CAPELLI \,IDENTITIES \,FOR \,LIE \,SUPERALGEBRAS}
\bigskip\bigskip
\centerline{\smc By Maxim NAZAROV}
\bigskip
\centerline{\hbox to 3.2cm{\hrulefill}}
\bigskip\bigskip


\centerline{\bf1.\ Introduction}\section{\,}
\kern-20pt

\nt
The Capelli identity [\hs1\hs] is one of the best exploited results            
of the classical invariant theory. It provides a set of
distinguished generators $C_1\lc C_N$ for the centre of the
enveloping algebra $\operatorname{U}(\frak{gl}_N)$
of the general linear Lie algebra. For any non-negative integer $M$
consider the natural action of the Lie algebra $\frak{gl}_N$
in the vector space $\CC^N\ot\CC^M$. Extend it to the action 
of the algebra $\operatorname{U}(\frak{gl}_N)$ in the
space of polynomial functions on $\CC^N\ot\CC^M$. The resulting
representation of $\operatorname{U}(\frak{gl}_N)$
by differential operators on $\CC^N\ot\CC^M$
with polynomial \text{coefficients is faithful when~$M\ge N$.} 

The image of the centre $\operatorname{Z}(\frak{gl}_N)$
of the algebra $\operatorname{U}(\frak{gl}_N)$ under this representation
coincides with the ring $\J$ of $\frak{gl}_N\x\frak{gl}_M$-invariant 
differential operators on $\CC^N\ot\CC^M$ with polynomial coefficients.
The latter ring has a distinguished set of generators $I_1\lc I_L$
with $L=\min\,(M,N)$ which are called the {\it Cayley operators} [\hs2\hs].
If $x_{ib}$ with $i=1\lc N$ and $b=1\lc M$ are the standard coordinates
on $\CC^N\ot\CC^M$ and $\d_{ib}$ are the corresponding partial
derivations, then $I_n$ equals
$$
\sum_{\,g\,\in\,S_n}\ \sum_{i_1\lc i_n}\sum_{b_1\lc b_n}\ 
{\operatorname{sgn}\,(g)}
\,\cdot\,
x_{i_{g(1)}b_1}\ldots x_{i_{g(n)}b_n} 
\,
\d_{i_1b_1}\ldots\,\d_{i_nb_n}\,.
\Tag{1.0}
$$
Here $\operatorname{sgn}\,(g)$ stands for the sign of the element $g$
of the symmetric group $S_n$. The Capelli identity gives
explicit formula for a preimage in $\operatorname{Z}(\frak{gl}_N)$
of the operator~$I_n$. Let $E_{ij}$ be the standard generators
of the enveloping algebra $\operatorname{U}(\frak{gl}_N)$ so that
in the above representation
$$
E_{ij}\,\mapsto\sum_{b}\,x_{ib}\,\d_{jb}\,.
$$
The Cayley operator $I_n$ is then the image of the element
$C_n\in\operatorname{Z}(\frak{gl}_N)$ equal to
$$
\sum_{\,g\,\in\,S_n}\ \sum_{i_1\lc i_n}\ 
{\operatorname{sgn}\,(g)}
\,\cdot\,
\prod^\rightarrow_s\ 
\bigl(\hskip1pt E_{i_{g(s)}i_s}+(s-1)\cdot\de_{i_{g(s)}i_s}\bigr)
\Tag{1.11}
$$
where the index $s$ runs through $1\lc n$ and factors in the
ordered product are arranged from the left to right while $s$
increases. Here $\de_{ij}$ is the Kronecker delta.

\vbox{
Let $\chi$ be the irreducible character
of the symmetric group $S_n$ corresponding to any partition of $n$
into not more than $M,N$ parts. Let us normalize $\chi$ so that
$\chi(1)=1$. When we replace the character
$\operatorname{sgn}$ in \(1.0) by $\chi$ we obtain an element
of distinguished basis in the vector space $\J\,$; see [\hs3\hs]. 
Explicit formula for a preimage in $\operatorname{Z}(\frak{gl}_N)$
of this basic element
was given in [\hs4\hs,\hs5\hs]. This result
generalizes classical formula \(1.11) and maybe called the
{\it higher Capelli identity}~[\hs4\hs]. In the present article
we extend this result to the queer Lie superalgebra $\qN$~[\hs6\hs].
In particular, we obtain an analogue for $\qN$ of the classical
Capelli identity. Such an analogue has been so far unknown.
See, however [\hskip1pt7\hskip1pt] and [\hskip1pt8\hskip1pt]
for related results.
}

\vbox{
The symmetric group $S_n$ appears in the formulas \(1.0) and \(1.11)
because its permutational action in the tensor product
$\End(\CC^N)^{\,\ot n}$
generates the commutant of the action of the Lie algebra $\frak{gl}_N$.
For the $n$-th tensor power of the defining representation
of the Lie superalgebra $\qN$ the role of $S_n$ is played [\hs9\hs] 
by the semidirect product $S_n\ltimes A_n\,$
where $A_n$ is the Clifford algebra with anticommuting
generators $a_1\lc a_n$. We will denote this
product by $\Sen$ and call it the Sergeev algebra. Note that
to obtain the higher Capelli identity one replaces $\operatorname{sgn}$ in
\(1.11) by a diagonal matrix element relative to the Young orthogonal
basis in the irreducible representation with character~$\chi$.
In Section 2 we give an analogue of this matrix element for
any irreducible representation of $\Sen$.

The irreducible $H_n$-modules
are parametrized by the partitions $\la$ of $n$ with
pairwise distinct parts. 
Note that the algebra $H_n$ has
a natural $\ZZ_2$-gradation, and we use the notion of a $\ZZ_2$-graded
irreducibility. Let $\ell_\la$ be the number of 
parts in $\la$. For each $\la$ we construct a certain element
$\Psi_\la\in\Sen$ such that under the left regular action of $\Sen$ the space
$\Sen\cdot\Psi_\la$ splits into a direct sum of $2^{[\ell_\la/2]}$ copies of
the irreducible module corresponding to $\la$; see Theorem 3.4 and
the subsequent remark. Moreover, $\Psi_\la$ is invariant with respect to
the natural involutive antiautomorphism of the algebra $\Sen$; see Lemma 2.3.
If $\la=(n)$ then

\vskip-12pt
$$
\Psi_\la\ =\prod_{1\le r<n}^\rightarrow
\biggl(\,\,\prod_{r<s\le n}^\rightarrow\ 
\biggl(
1-\frac{(r\ts s)}{u_r-u_s}+\frac{(r\ts s)\cdot a_r\ts a_s}{u_r+u_s}
\biggr)\biggr)
\Tag{1.0101}
$$
\vskip-4pt
\nt
where $u_s=\sqrt{s(s-1)}$
while $(r\hskip1pt s)\in S_n$ is the transposition of $r$ and $s$.
To give an explicit formula for $\Psi_\la$ with a general $\la$
we use the fusion procedure [\hs10\hs]; see Theorem 2.2 here.
More detailed exposition of this construction appears in
\text{[\hs11\hs].}

In Section 3 we introduce our main technical tool -- the Jucys-Murphy
elements of the algebra $\Sen$. These are the pairwise commuting elements
$x_1\lc x_n$ defined by (3.1); see also Lemma 3.2.
With respect to the left regular action of $H_n$ the element $\Psi_\la$
is a joint eigenvector of $x_1\lc x_n$ and the corresponding
eigenvalues are easy to describe; see Proposition 3.3. Our proof
of the Capelli identity for 
$\qN$ will be
based on Proposition 3.6\hskip1pt;
cf. [\hs12\hs,\hs13\hs]. Namely, we will use Corollary 3.7.

Section 4 contains the main result of this article. We define
actions of the Lie superalgebras $\qN$ and $\qM$ in the free
supercommutative algebra $\PMN$
with the even generators $x_{ib}$
and odd generators $x_{-i,b}$ where $i=1\lc N$ and $b=1\lc M$;
see Proposition 4.1. Let $\PDMN$ be the algebra generated
by the operators of left multiplication in $\PMN$ by $x_{\pm i,b}$
and by the corresponding left derivations $\d_{\pm i,b}$.
So we get a representation $\ga:\UN\to\PDMN$ for the enveloping
algebra of $\qN$. The image of the centre $\ZN\subset\UN$ with respect to
$\ga$ coincides with the ring $\I$ of $\qN\x\qM$-invariants in $\PDMN$.
We introduce a distinguished basis in the vector space $\I$;
see Proposition 4.3. The elements of this basis are parametrized
by the partitions $\la$ of $n=0,1,2,\ldots$ with $\ell_\la\le M,N$
and are determined by (4.7). Our main result is the explicit formula (4.6)
for a preimage $C_\la\in\ZN$ of the basic element $I_\la\in\I$
corresponding to $\la$; see also (4.5). 
The equality $I_\la=\ga\,(C_\la)$ with $\la=(n)$
maybe regarded as an analogue for $\qN$ of the classical Capelli identity.

I am very grateful to I.\,Cherednik for many illuminating conversations.
I am also {grateful} to I. Penkov for his interest in this work.
I am especially indebted to G.\,Olshanski.
Discussion with him of the results [\hs4\hs] has inspired
\text{the present work.}
}


\newpage

\centerline{\bf2.\ Fusion procedure for the Sergeev algebra}\section{\,}
\kern-20pt

\nt
We start with recalling several known facts about irreducible
modules over the {\it Sergeev algebra} $\Sen$.
By definition, $\Sen$ is the semidirect product
of the symmetric group $S_n$ and the Clifford algebra $A_n$
with $n$ generators over the complex field $\CC$. 
These generators are denoted by $a_1\lc a_n$ and subjected to the relations
$$
a_i^2=-1\,;
\qquad
a_i\,a_j=-a_j\,a_i\,,
\quad
i\neq j.
$$
The symmetric group $S_n$ acts on the algebra $A_n$ by permutations
of these $n$ generators. Denote by the superscript $\ast$
the involutive
antiautomorphism of the algebra $\Sen$ defined by the assignments
$a_i\mapsto a_i\1$ and $g\mapsto g\1$ for any $g\in S_n$.

For any $\ZZ_2$-graded algebra $\operatorname{A}$
a representation $\operatorname{A}\to\End(\CC^{K|K})$
will be called {\it irreducible} 
if the even part of its supercommutant equals $\CC$.
If the 
supercommutant coincides with $\CC$ this representation is called 
{\it absolutely irreducible}. 
We will equip the algebra $\Sen$ with $\ZZ_2$-gradation so that
$\deg a_i=1$ and $\deg g=0$ for any element $g\in S_n\,$.
The irreducible modules over
the $\ZZ_2$-graded algebra $\Sen$ are parametrized by partitions 
$\la$ of $n$ with pairwise distinct parts.
Such a partition is called {\it strict\,}. The $\Sen$-module $U_\la$
corresponding to $\la$ is absolutely irreducible if and only if
the number $\ell_\la$ of non-zero parts in $\la$ is even
[\hs9\hs,\hs Lemma 6\hs].

Consider the left regular representation of the algebra $\Sen$.
In this section for any strict partition $\lambda$ we will construct
a certain element $\Psi_\la\in \Sen$ such that the left ideal
$\Sen\cdot\Psi_\la$
is a direct sum of $2^{[\ell_\la/2]}$ copies of the irreducible 
$\Sen$-module corresponding to $\la$. Moreover, we will have
the equality $\Psi_\la^\ast=\Psi_\la$. Our construction is motivated by
the results of [\hs10\hs] and [\hs14\hs]; see [\hs11\hs] for more detail.

Strict partitions are usually depicted as {shifted Young diagrams}.
For instance, here is the diagram corresponding to the partition
$\lambda=(4,3,1)$:

\smallskip\smallskip
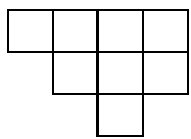
\vsk>
\vbox{
$$
{\bx}
{\bx}
{\bx}
{\bx}
$$
\vglue-17.8pt
$$
\phantom{\bx}
{\bx}
{\bx}
{\bx}
$$
\vglue-17.8pt
$$
\phantom{\bx}
\phantom{\bx}
{\bx}
\phantom{\bx}
$$
}
\smallskip\smallskip\smallskip
\nt
We will denote by $\La$ the shifted {\it column tableau} of the shape $\la$.
It is obtained by filling the boxes of $\la$ with the numbers $1\lc n$
by columns from the left to the right, downwards in every column.
For each $i=1\lc n$ we put $c_i=p-q$ if the number $i$ appears in the
$p$-th column and $q$-th row of the tableau $\La$. The difference
$p-q$ is then called the {\it content} of the box of the diagram $\la$
occupied by the number $i$. For example, here on the left we show
the shifted column tableau of the shape $\la=(4,3,1)$:

\smallskip\smallskip
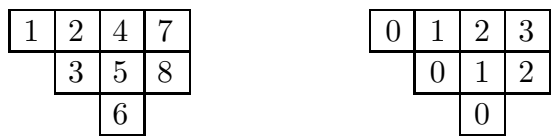
\vsk>
\vbox{
$$
{\bx}
{\bx}
{\bx}
{\bx}
\phantom{\bx}
\phantom{\bx}
\phantom{\bx}
\phantom{\bx}
{\bx}
{\bx}
{\bx}
{\bx}
$$
\vglue-17.8pt
$$
\phantom{\bx}
{\bx}
{\bx}
{\bx}
\phantom{\bx}
\phantom{\bx}
\phantom{\bx}
\phantom{\bx}
\phantom{\bx}
{\bx}
{\bx}
{\bx}
$$
\vglue-17.8pt
$$
\phantom{\bx}
\phantom{\bx}
{\bx}
\phantom{\bx}
\phantom{\bx}
\phantom{\bx}
\phantom{\bx}
\phantom{\bx}
\phantom{\bx}
\phantom{\bx}
{\bx}
\phantom{\bx}
$$
\vglue-57.5pt
$$
1
\kern9pt
2
\kern9pt
4
\kern9pt
7
\kern67pt
0
\kern9pt
1
\kern9pt
2
\kern9pt
3
$$
\vglue-18pt
$$
\phantom1
\kern9pt
3
\kern9pt
5
\kern9pt
8
\kern67pt
\phantom0
\kern9pt
0
\kern9pt
1
\kern9pt
2
$$
\vglue-18pt
$$
\phantom1
\kern9pt
\phantom3
\kern9pt
6
\kern9pt
\phantom8
\kern67pt
\phantom0
\kern9pt
\phantom0
\kern9pt
0
\kern9pt
\phantom0
$$
\smallskip\smallskip\smallskip\smallskip
\nt
On the right we indicated the contents of the boxes of the shifted 
Young diagram.
}

For any distinct $i,j=1\lc n$ let $(i\ts j)$ be the transposition
in the symmetric group $S_n$. Consider the rational function of two complex
variables $u,v$
valued in the algebra $\Sen$
$$
\phi_{ij}(u,v)=1-\frac{(i\ts j)}{u-v}+\frac{(i\ts j)\cdot a_i\ts a_j}{u+v}.
$$
As a direct calculation shows, this rational function satisfies the
equations
$$
\phi_{ij}(u,v)\ts\ts\phi_{ik}(u,w)\ts\ts\phi_{jk}(v,w)=
\phi_{jk}(v,w)\ts\ts\phi_{ik}(u,w)\ts\ts\phi_{ij}(u,v) 
\Tag{2.1}
$$
for all pairwise distinct $i,j,k$.
Evidently, for all pairwise distinct $i,j,k,l$ we have
$$
\phi_{ij}(u,v)\ts\ts\phi_{kl}(z,w)=\phi_{kl}(z,w)\ts\ts\phi_{ij}(u,v)
\Tag{2.2}
$$
We will also make use of the relations for all distinct $i,j$
$$
\gather
\phi_{ij}(u,v)\,\phi_{ji}(v,u)=1-\frac1{(u-v)^2}-\frac1{(u+v)^2}\,;
\Tag{2.31}
\\
a_i\,\phi_{ij}(u,v)\,a_i\1=\phi_{ij}(-\,u,v)\,,
\Tag{2.32}
\\
a_j\,\phi_{ij}(u,v)\,a_j\1=\phi_{ij}(u,-\,v)\,.
\Tag{2.33}
\endgather
$$
Note that due to \(2.31) the element $\phi_{ij}(u,v)\in \Sen$ is not
invertible
if and only if 
$$
\frac1{(u-v)^2}+\frac1{(u+v)^2}=1\,.
\Tag{2.34}
$$
Observe also that in the latter case the element $\phi_{ij}(u,v)/2\in \Sen$
is an idempotent.

Consider the rational function of $u,v,w$ appearing at either side
of \(2.1). Denote by $\phik$ this function.
The factor $\phi_{ik}(u,w)$ in \(2.1) has a pole at $u=\pm w$.
Still we have the following lemma. 
It will be the basis for  constructing $\Psi_\la\in \Sen$.

\proclaim{Lemma 2.1}
Restriction of $\phik$ to the set of $(u,v,w)$
such that the pair $(u,v)$ satisfies the condition \(2.34),
is continuous at $u=\pm w$.
\endproclaim

\demo{Proof}
Assume that the the condition \(2.34) is satisfied.
Then due to \(2.31) and \(2.33)
$$
\phi_{ij}(u,v)\,\phi_{ji}(v,u)=0\,,
\qquad
\phi_{ij}(u,v)\,a_i\,\phi_{ji}(v,-u)=0\,.
$$
Hence the product $\phi_{ijk}(u,v,w)$ can be rewritten as
$$
\gather
\phi_{ij}(u,v)\,\phi_{jk}(v,w)
-\phi_{ij}(u,v)\,\frac{(ik)}{u-w}\,
\bigl(\phi_{jk}(v,w)-\phi_{jk}(v,u)\bigr)\,+
\\
\phi_{ij}(u,v)\,\frac{(ik)\,a_i\,a_k}{u+w}\,
\bigl(\phi_{jk}(v,w)-\phi_{jk}(v,-u)\bigr)\,.
\endgather
$$
Under the condition \(2.34) the latter  function is continuous
at $u=\pm w$
\enddemos

\nt
For any two real non-negative variables $s,t$ let us
substitute in the equation \(2.34)
$$
u=\sqrt{s(s+1)}\,,
\quad
v=\sqrt{t(t+1)}\,.
$$
Observe that \(2.34) will be then satisfied if $s-t=\pm1$. We will
denote for short
$$
\gather
\psi_{ij}(s,t)=\phi_{ij}\bigl(\sqrt{s(s+1)},\sqrt{t(t+1)}\,\bigr)\,,
\\
\psi_{ijk}(s,t,r)=\phi_{ijk}
\bigl(\sqrt{s(s+1)},\sqrt{t(t+1)},\sqrt{r(r+1)}\,\bigr)\,.
\endgather
$$ 

Now introduce a real non-negative parameter $t_i$ for each $i=1\lc n$.
Equip the set of all pairs $(i,j)$ where $1\le i<j\le n$
with the lexicographical ordering. Introduce the ordered
product over this set
$$
\prod_{(i,j)}^\rightarrow\ts\ts
\psi_{ij\ts}(c_i+t_i\ts,c_j+t_j).
\Tag{2.3}
$$
Consider this product as a function of the parameters
$t_1\lc t_n$ valued in the algebra $\Sen$.
Le us denote by $\Psiz$ this function.
It may have singularities when $c_i+t_i=c_j+t_j$ for some $i\neq j$.
Consider 
the set $\Z$ of all tuples $(t_1\lc t_n)$
such that $t_i=t_j$ whenever the numbers
$i$~and~$j$ appear in the same row of the tableau $\La$.
So $\Z=\Bbb R_{\geqslant0}^{\ell_\la}$.
The following
theorem goes back \text{to [\hs10\hs] and [\hs11\hs].}

\proclaim{Theorem 2.2}
Restriction of $\Psiz$ to $\Z$ is continuous at
$t_1=\ldots=t_n$. 
\endproclaim

\demo{Proof}
We shall provide an expression for
the restriction of the function \(2.3) to $\Z$ which is manifestly continuous
at $t_1=\ldots=t_n$.
Let us reorder the pairs $(i,j)$ in the product \(2.3) as follows.
This reordering will not affect the value of the product due to
the relations \(2.1) and \(2.2).
Let $\C$ be the sequence of numbers obtained by reading the tableau
$\La$ in the usual way, that is by rows from the top to the bottom, eastwards
in every row. For~each $j=1\lc n$ denote by $\A_j$ and $\B_j$
the subsequences of $\C$ consisting of all numbers $i<j$ which appear
respectively before and after $j$ in that sequence.
Now set $(i,j)\prec(k,l)$ if one of the following conditions is satisfied:

\medskip
\vbox{
\line{\ -\ the number $i$ appears in $\A_j$ while $k$ appears in $\B_l$;
\hfill}
\line{\ -\ the numbers $i$ and $k$ appear respectively in $\A_j$ and $\A_l$
where $j<l$;
\hfill}
\line{\ -\ the numbers $i$ and $k$ appear respectively in $\B_j$ and $\B_l$
where $j>l$;
\hfill}
\line{\ -\ we have the equality $j=l$ and $i$ appears before $k$ in $A_j$
or $B_j$.
\hfill}
}
\smallskip

{}From now on we assume that the factors in \(2.3) corresponding to the
pairs $(i,j)$ are arranged with respect to this new ordering.
The factor $\psi_{\ts i\ts j\ts}(c_i+t_i\ts,c_j+t_j)$ has a singularity
at $t_i=t_j$
if and only if  $i$ and $j$ stand on the same diagonal of the
tableau $\La$. We will then call the pair $(i,j)$ {\it singular.}
Observe that the number $i$ occurs in the subsequence $\B_j$exactly when
$i$ stands
to the left and below of $j$ in the tableau $\La$.
In this case $c_j-c_i>1$ and the pair $(i,j)$ cannot be singular. 

Let a singular pair $(i,j)$ be fixed. Suppose that the number $i$ appears
in the $p$-th column and the $q$-th row of the tableau $\La$.
In our new ordering the next pair after $(i,j)$ is $(h,j)$ where the number
$h$
appears in the $(p+1)$-th column and the $q$-th row of
$\La$. In particular, we have $c_i=c_j=c_h-1$. Moreover, $(i,h)\prec(i,j)$.
Due to the relations \(2.1)\ts,\(2.2) the product
$$
\prod_{(k,l)\prec(i,j)}^\rightarrow\ts
\psi_{kl\ts}(c_k+t_k\ts,c_l+t_l)
$$
is divisible on the right by
$\psi_{ih\ts}(c_i+t_i\ts,c_h+t_h)$.
Note that each value of the restriction of 
$\psi_{ih\ts}(c_i+t_i\ts,c_h+t_h)/2$
to $t_i=t_h$ is an idempotent in $\Sen$.

Now for each singular pair $(i,j)$ let us replace the two adjacent factors
in \(2.3)
$$
\psi_{ij\ts}(c_i+t_i\ts,c_j+t_j)
\,
\psi_{hj\ts}(c_h+t_h\ts,c_j+t_j)
$$
by  
$$
\gather
\psi_{ih\ts}(c_i+t_i\ts,c_h+t_h)
\,
\psi_{ij\ts}(c_i+t_i\ts,c_j+t_j)
\,
\psi_{hj\ts}(c_h+t_h\ts,c_j+t_j)/2=
\\
\psi_{ihj\ts}(c_i+t_i\ts,c_h+t_h\ts,c_h+t_h)/2\ts.
\Tag{2.35}
\endgather
$$
This replacement does not affect the value of restriction to $\Z$ of the
function \(2.3). But the restriction to $t_i=t_h$ of  \(2.35)
is continuous at $t_i=t_j$ by Lemma 2.1
\enddemos

\nt
The process of continuation of the function $\Psiz$ along the set $\Z$
is called the {\it fusion procedure}.
In the proof of Theorem 2.2 we established the decomposition
$$
\Psiz=\Upz\cdot\Thez
$$
where $\Upz$ and $\Thez$ are products of the factors in \(2.3) which
correspond
to the pairs $(i,j)$ with $i$ appearing in $\A_j$ and $\B_j$ respectively.
The function $\Thez$ is continuous at $t_1=\ldots=t_n$. Moreover, any value
of this function at $t_1=\ldots=t_n$ is invertible.
Let us denote by $\The_\la$ the value $\The_\la(0\lc 0)$.
Restriction to $\Z$ of the function $\Upz$ is continuous at $t_1=\ldots=t_n$ 
as well as the restriction of the function $\Psiz$. Denote
respectively by $\Up_\la$ and $\Psi_\la$ the values of these restrictions
at $t_1=\ldots=t_n=0$. Then $\Psi_\la=\Up_\la\ts\The_\la$.
Let $\alpha$ be the linear map $\Sen\to A_n$
identical on $A_n$ such that
$\al(ga)=0$ for  $g\neq1$ in $S_n$ and any element $a\in A_n$.

\proclaim{Lemma 2.3}
We have $\Psi_\la^\ast=\Psi_\la$ and $\al(\Psi_\la)=1$.
\endproclaim

\demo{Proof}
By the definition of the antiautomorphism $\ast$ we have 
$\phi_{ij}(u,v)^\ast=\phi_{ij}(u,v)$ for any distinct indices $i$ and $j$.
Therefore due to the relations \(2.1) and \(2.2)
the product~\(2.3) is invariant with respect to this
antiautomorphism. So is the value $\Psi_\la$ of its restriction to $\Z$.
Further, we have the equality $\al\bigl(\Psi_\la(t_1\lc t_n)\bigr)=1$
by the definition \(2.3). Hence $\al(\Psi_\la)=1$ 
\enddemos

\nt
Throughout this article we will denote $z_i=\sqrt{c_i(c_i+1)}$
for $i=1\lc n$.

\proclaim{Proposition 2.4}
Let the numbers $k<l$ stand next to each other in one row of the tableau
$\La$.
Then the element $\Up_\la\in \Sen$ is divisible on the right by
$\phi_{kl}(z_k,z_l)$.
\endproclaim

\demo{Proof}
Due to the relations \(2.1) and \(2.2) the restriction of $\Upz$ to $\Z$ is
divisible
on the right by $\psi_{kl\ts}(c_k+t_k\ts,c_l+t_l)$. Here $t_k=t_l$ 
and the element $\psi_{kl\ts}(c_k+t_k\ts,c_l+t_l)/2\in \Sen$ is an
idempotent.
Restriction of $\Upz$ to $\Z$ is continuous at
$t_1=\ldots=t_n=0$. So $\Up_\la$ is divisible on right by
$\psi_{kl}(c_k,c_l)/2$
\enddemos

\proclaim{Corollary 2.5}
Let the numbers $k<l$ stand next to each other in the first row of the
tableau
$\La$. Then the element $\Psi_\la\in \Sen$ is divisible on the right by
$$
\phi_{\ts kl\ts}(z_k\ts,z_l)\ts\ts\cdot\hskip-2pt
\prod_{k<m<l}^\rightarrow\ts
\phi_{\ts ml\ts}(z_m\ts,z_l).
\Tag{2.5}
$$
The element $\Psi_\la$ is then also divisible on the left by
$$
\prod_{k<m<l}^\leftarrow\ts
\phi_{\ts ml\ts}(z_m\ts,z_l)
\cdot
\phi_{\ts kl\ts}(z_k\ts,z_l).
\Tag{2.6}
$$
\endproclaim

\demo{Proof}
By Proposition 2.4 the element $\Up_\la$ is divisible on the right by
$\psi_{\ts kl\ts}(c_k\ts,c_l).$ But due to the relations \(2.1) and \(2.2)
the product $\psi_{\ts kl\ts}(c_k\ts,c_l)\cdot\The_\la$ is divisible on the
right
by \(2.5). Since $\Psi_\la=\Up_\la\The_\la$ we obtain the first statement
of Corollary~2.5. Note that the image of \(2.5) with respect to the
antiautomorphism
$\ast$ is \(2.6). Since by Lemma 2.3 the element $\Psi_\la$
is invariant with respect to this antiautomorphism
we get the second statement of Corollary 2.5
\enddemos

\nt
Consider again the rational function $\phik$ appearing at either side of
\(2.1).
The values at $u=w$ of its restriction to $(u,v)$ subjected to \(2.34)
need not be divisible on the right by $\phi_{jk}(v,u)=\phi_{jk}(v,w)$.
Yet we have the following lemma.

\proclaim{Lemma 2.6}
Restriction of the function $\phi_{ijk}(u,v,w)\,\phi_{kj}(w,v)$
to those $(u,v)$ which satisfy \(2.34), takes at $u=w$ the value
$$
(ik)\,\phi_{kj}(w,v)\cdot\biggl(\frac2{(v+w)^3}-\frac2{(v-w)^3}\biggr)\,.
$$
\endproclaim

\demo{Proof}
It consists of a direct calculation. Namely, by our definition
$$
\phi_{ijk}(u,v,w)\,\phi_{kj}(w,v)=\phi_{ij}(u,v)\,\phi_{ik}(u,w)\cdot
\phi_{jk}(v,w)\,\phi_{kj}(w,v).
$$
Under the condition \(2.34) we have by \(2.31) the equality
$$
\phi_{jk}(v,w)\,\phi_{kj}(w,v)=
\frac1{(u\!-\!v)^2}+\frac1{(u\!+\!v)^2}-\frac1{(v\!-\!w)^2}
-\frac1{(v\!+\!w)^2}\,.
\Tag{2.61}
$$
Dividing the right hand side of \(2.61) by $u^2-w^2$ and then
setting $u=w$ we get
$$
-\,\frac1w\cdot\biggl(\frac1{(v+w)^3}-\frac1{(v-w)^3}\biggr).
$$
On the other hand, by setting $u=w$ in the product 
$$
\align
&\phi_{ij}(u,v)\,\phi_{ik}(u,w)\cdot(u^2-w^2)=
\\
&\phi_{ij}(u,v)\cdot\bigl((u^2-w^2)-(ik)\,(u+w)+(ik)\,a_ia_k\,(u-w)\bigr)
\endalign
$$
we obtain
$$
-\,2\,w\cdot\phi_{ij}(u,v)\,(ik)=-\,2\,w\cdot(ik)\,\phi_{kj}(w,v)
\quad\square
$$
\enddemo

\proclaim{Corollary 2.7}
Restriction of the function $\phi_{ijk}(u,v,w)\,\phi_{kj}(w,v)$
to those $(u,v)$ which satisfy \(2.34), vanishes at $u=w=0$.
\endproclaim

\nt
The next proposition makes the central part of the present section;
cf. [\hs15\hs].

\proclaim{Proposition 2.8}
Let the numbers $k$ and $k+1$ stand in the same column of the tableau $\La$.
Then the element $\Up_\la\in \Sen$ is divisible on the left by
$\phi_{k,k+1}(z_k,z_{k+1})$.
\endproclaim

\demo{Proof}
Observe first that Proposition 2.8 follows from its particular case $k+1=n$.
Indeed, let $\nu$ be the shape of the tableau obtained from $\La$ by removing
each of the numbers $k+2\lc n$. Then
$$
\Upz=\Up_\nu(t_1\lc t_{k+1})
\cdot
\prod_{(i,j)}^\rightarrow\ts\ts
\psi_{ij\ts}(c_i+t_i\ts,c_j+t_j)
$$
where $j=k+2\lc n$ and $i$ runs through the sequence $\A_j$.
Consider the value $\Up_\nu$ at $t_1=\ldots=t_{k+1}=0$
of the restriction of $\Up_\nu(t_1\lc t_{k+1})$ to $\Z$.
According to our proof of Theorem 2.2 then
$\Up_\la=\Up_\nu\ts\Up$ for a certain element $\Up\in \Sen$.

{}From now on we will assume that $k+1=n$. Since the element $\The_\la$
is invertible, it suffices to prove that $\Up_\la\ts\The_\la=\Psi_\la$
is divisible by $\psi_{n-1,n}(c_{n-1},c_n)$ on the left. But
$\Psi_\la^\ast=\Psi_\la$ by Lemma 2.3. So we will prove that  $\Psi_\la$
is divisible by $\psi_{n-1,n}(c_{n-1},c_n)$ on the right.
Since $\psi_{n-1,n}(c_{n-1},c_n)+\psi_{n\ts,n-1\ts}(c_n\ts,c_{n-1})=2$,
this is to prove
$$
\Psi_\la\cdot\psi_{n\ts,n-1\ts}(c_n\ts,c_{n-1})=0.
\Tag{2.7}
$$

Suppose that the number $n$ appears in the $p$-th column and the $q$-th row
of the tableau $\La$. Then by our assumption the number $n-1$ appears in the
same column and $(q-1)$-th row of $\La$. 
Let $i_1<\ldots<i_r$ be all the numbers in the $q$-th row.
So we have $i_r=n$.
Then due to the relations \(2.1) and \(2.2) we have for a certain element
$\The\in \Sen$ the equality
$$
\The_\la\cdot\psi_{n\ts,n-1\ts}(c_n\ts,c_{n-1})=
\prod_{s<r}^\rightarrow\ts
\psi_{\ts i_s,\ts n-1\ts}(c_{i_s}\ts,c_{n-1}) \cdot
\psi_{n\ts,n-1\ts}(c_n\ts,c_{n-1})
\cdot
\The.
$$
Therefore to get \(2.7)
we have to prove that
$$
\Up_\la\ts\ts\cdot\ts
\prod_{s<r}^\rightarrow\ts
\psi_{\ts i_s,\ts n-1\ts}(c_{i_s}\ts,c_{n-1}) \cdot
\psi_{n\ts,n-1\ts}(c_n\ts,c_{n-1})=0.
\Tag{2.8}
$$

We will now prove \(2.8) by induction on $r$. Suppose that $r=1$.
Let $m$ be the number appearing in the $(p-1)$-th column and $(q-1)$-th row
of $\La$.
Then according to our proof of Theorem 2.2 the function
$\Upz$ has~the~form
$$
\Up(t_1\lc t_n)\cdot
\psi_{mn}(c_m+t_m,c_n+t_n)\,
\psi_{n-1,n}(c_{n-1}+t_{n-1},c_n+t_n)
$$
where the restriction of $\Up(t_1\lc t_n)$ to $\Z$ is continuous
at $t_1=\ldots=t_n=0$.
Moreover, this restriction is divisible on the right by 
$\psi_{m,n-1}(c_m+t_m,c_{n-1}+t_{n-1})$ where $t_m=t_{n-1}$.
Since $c_m=c_n=0$ and $c_{n-1}=1$, restriction to $t_m=t_{n-1}$ of 
$$
\psi_{m,n-1,n}(c_m+t_m,c_{n-1}+t_{n-1},c_n+t_n)\cdot
\psi_{n\ts,n-1\ts}(c_n+t_n\ts,c_{n-1}+t_{n-1})
$$
vanishes at $t_m=t_n=0$ by Corollary 2.7. This proves \(2.8)
for $r=1$.

Now suppose that $r>1$. We have to prove that the restriction to $\Z$ of
$$
\Upz\ts\ts\cdot\ts
\prod_{s\le r}^\rightarrow\ts
\psi_{\ts i_s,\ts n-1\ts}(c_{i_s}+t_{i_s}\ts,c_{n-1}+t_{n-1})
\Tag{2.10}
$$
vanishes at $t_1=\ldots=t_n=0$. Now denote $i_{r-1}=m$. The number $m-1$
appears in the $(p-1)$-\ts th column and the $(q-1)$-\ts th row of $\La$.
So we have $c_{m-1}=c_n$.
Let $\mu$ be the shape of the tableau obtained from $\La$ by removing each
of the numbers $m+1\lc n$. Then the function $\Upz$ has the form
$$
\gather
\Up_\mu(t_1\lc t_m)\cdot\Psi(t_1\lc t_{n-1}) \cdot
\psi_{\ts m-1,\ts n-1\ts}(c_{m-1}+t_{m-1}\ts,c_{n-1}+t_{n-1})\times
\\
\Phi(z_1\lc z_{n})
\cdot
\psi_{\ts m-1,\ts n\ts}(c_{m-1}+t_{m-1}\ts,c_{n}+t_{n}) \,
\psi_{\ts n-1,\ts n\ts}(c_{n-1}+t_{n-1}\ts,c_{n}+t_{n})\times
\\
\prod_{s<r}^\rightarrow\ts
\psi_{\ts i_sn\ts}(c_{i_s}+t_{i_s}\ts,c_{n}+t_{n}).
\endgather
$$
Here we have denoted by $\Psi(t_1\lc t_{n-1})$ the product
$$
\prod_{(i,j)}^\rightarrow\ts\ts
\psi_{ij\ts}(c_i+t_i\ts,c_j+t_j)\ts;
\qquad
j=m+1\lc n-1
\Tag{2.85}
$$
where $i$ runs through $\A_j$ but
$(i,j)\neq(m-1\ts,n-1)$. Further, we have denoted
$$
\Phi(t_1\lc t_{n})=
\prod_{(i,n)}^\rightarrow\ts\ts
\psi_{in\ts}(c_i+t_i\ts,c_n+t_n)
\Tag{2.9}
$$
where $i$ runs through the sequence $\A_n$ but $i\neq m-1\ts,n-1\ts\lc m$.
In particular, any factor in the product \(2.9) commutes with
$$
\psi_{\ts m-1,\ts n-1\ts}(c_{m-1}+t_{m-1}\ts,c_{n-1}+t_{n-1})
$$
due to \(2.2). Therefore the product \(2.10) takes the form
$$
\gather
\Up_\mu(t_1\lc t_m)
\cdot
\Psi(t_1\lc t_{n-1})
\cdot
\Phi(t_1\lc t_{n})
\,\times
\Tag{2.11}
\\
\psi_{\ts m-1,\ts n-1\ts,n\ts}
(c_{m-1}+t_{m-1}\ts,c_{n-1}+t_{n-1},c_{n}+t_{n})
\,\times
\\
\prod_{s<r}^\rightarrow\ts
\psi_{\ts i_sn\ts}(c_{i_s}+t_{i_s}\ts,c_{n}+t_{n}) \cdot
\prod_{s<r}^\rightarrow\ts
\psi_{\ts i_s,\ts n-1\ts}(c_{i_s}+t_{i_s}\ts,c_{n-1}+t_{n-1})
\,\times
\\
\psi_{n\ts,n-1\ts}(c_n+t_n\ts,c_{n-1}+t_{n-1})=
\\
\Up_\mu(t_1\lc t_m)
\cdot
\Psi(t_1\lc t_{n-1})
\cdot
\Phi(t_1\lc t_{n})
\,\times
\\
\psi_{\ts m-1,\ts n-1\ts,n\ts}(c_{m-1}+t_{m-1}\ts,c_{n-1}
+t_{n-1},c_{n}+t_{n})
\,\,
\psi_{n\ts,n-1\ts}(c_n+t_n\ts,c_{n-1}+t_{n-1})
\\
\prod_{s<r}^\rightarrow\ts
\psi_{\ts i_s,\ts n-1\ts}(c_{i_s}+t_{i_s}\ts,c_{n-1}+t_{n-1})
\cdot\prod_{s<r}^\rightarrow\ts
\psi_{\ts i_sn\ts}(c_{i_s}+t_{i_s}\ts,c_{n}+t_{n})\,.
\endgather
$$
To get the latter equality we used the relations \(2.1) and \(2.2).
Restriction to $\Z$ of the product of factors in the first line of \(2.11)
is continuous at $t_1=\ldots=t_n=0$ according to our proof of Theorem 2.2.
Each of the factors in the last line is also continuous at
$t_1=\ldots=t_n=0$.
Therefore by Lemma 2.6 the restriction of \(2.11) to $\Z$ has 
at $t_1=\ldots=t_n=0$ the same value
as restriction to $\Z$ of the product
$$
\gather
\Up_\mu(t_1\lc t_m)
\cdot
\Psi(t_1\lc t_{n-1})
\cdot
\Phi(t_1\lc t_{n})
\,\times
\Tag{2.12}
\\
(m-1,n)\cdot
\phi_{n\ts,n-1\ts}(c_n+t_n\ts,c_{n-1}+t_{n-1})\cdot f\,\,
\,\times
\\
\prod_{s<r}^\rightarrow\ts
\phi_{\ts i_s,\ts n-1\ts}(c_{i_s}+t_{i_s}\ts,c_{n-1}+t_{n-1})
\cdot\prod_{s<r}^\rightarrow\ts
\phi_{\ts i_s,\ts n\ts}(c_{i_s}+t_{i_s}\ts,c_{n}+t_{n})=
\\
\Up_\mu(t_1\lc t_m)
\cdot
\Psi(t_1\lc t_{n-1})
\cdot
\Phi(t_1\lc t_{n})
\,\times
\\
\prod_{s<r}^\rightarrow\ts
\psi_{\ts i_s,\ts m-1\ts}(c_{i_s}+t_{i_s}\ts,c_{n}+t_{n})
\cdot
\prod_{s<r}^\rightarrow\ts
\psi_{\ts i_s,\ts n-1\ts}(c_{i_s}+t_{i_s}\ts,c_{n-1}+t_{n-1})
\,\times
\\
(m-1,n)\cdot \psi_{n,n-1\ts}(c_n+t_n\ts,c_{n-1}+t_{n-1})\cdot f
\endgather
$$
for a certain number $f\in\Bbb R$. Here each of the factors
$\psi_{\ts i_s,\ts m-1\ts}(c_{i_s}+t_{i_s}\ts,c_{n}+t_{n})$ commutes with
$\Phi(t_1\lc t_{n})$ by the relations \(2.2).
In each of these factors we can replace $c_{n}+t_{n}$ by $c_{m-1}+t_{m-1}$
without affecting the value at $t_1=\ldots=t_n=0$ of the restriction to $\Z$ 
of \(2.12). Denote
$$
\Gaz\ts=\ts
\prod_{s<r}^\rightarrow\ts
\psi_{\ts i_s,\ts m-1\ts}(c_{i_s}+t_{i_s}\ts,c_{m-1}+t_{m-1}).
$$
According to our proof of Theorem 2.2 it now suffices to demonstrate
vanishing at $t_1=\ldots=t_{n-1}$ of the restriction to $\Z$ of the product
$$
\Up_\mu(t_1\lc t_m)
\cdot
\Psi(t_1\lc t_{n-1})
\cdot
\Gaz
\ts.
\Tag{2.13}
$$

Consider the product \(2.85). Here the factors corresponding to the pairs
$(i,j)$ are arranged with respect to ordering chosen in the proof of
Theorem 2.2.
Let us now reorder the pairs $(i,j)$ in \(2.85) as follows. For each
number $j>m$
appearing in the $(p-1)$-th column of $\La$ change the sequence
$$
(m-1,j)\ts,(i_1,j)\lc(i_{\ts r-1},j)
$$
to
$$
(i_1,j)\lc(i_{\ts r-1},j)\ts,(m-1,j)\ts.
$$
Denote by $\Psi^{\,\prime}(t_1\lc t_{n-1})$ the resulting ordered product.
Then by \(2.1) and \(2.2)
$$
\Psi(t_1\lc t_{n-1})\cdot\Gaz=
\Gaz\cdot\Psi^{\,\prime}(t_1\lc t_{n-1}).
\Tag{2.135}
$$

Now let $(i,j)$ be any singular pair in \(2.85).
Let $(h,j)$ be the pair following $(i,j)$ in the ordering from our proof of
Theorem 2.2. Then $(h,j)$ follows $(i,j)$ in our new ordering as well.
Furthermore, by the relations \(2.1) and \(2.2) the product
$$
\Up_\mu(t_1\lc t_m)
\cdot
\Gaz
\Tag{2.14}
$$
is divisible on the right by $\psi_{ih}(c_i+t_i,\ts c_h+t_h)$. Therefore
$$
\align
&
\Up_\mu(t_1\lc t_m)
\cdot
\Gaz
\cdot
\Psi^{\,\prime}(t_1\lc t_{n-1})=
\Tag{2.15}
\\
&
\Up_\mu(t_1\lc t_m)
\cdot
\Gaz
\cdot
\Psi^{\,\prime\prime}(t_1\lc t_{n-1})
\endalign
$$
where restriction of $\Psi^{\,\prime\prime}(t_1\lc t_{n-1})$ to $\Z$
is continuous at $t_1=\ldots=t_{n-1}=0$.
But by the inductive assumption the restriction of the product
\(2.14) to $\Z$ vanishes
at $t_1=\ldots=t_{m}=0$.
Thus due to \(2.135) and \(2.15) the restriction of \(2.13) to $\Z$ vanishes
at $t_1=\ldots=t_{n-1}=0$ as well
\enddemos

\nt
We have shown that the element $\Psi_\la=\Up_\la\,\The_\la$ of $\Sen$
is non-zero and $\ast$-invariant. 

\proclaim{Corollary 2.9}
Let the numbers $k$ and $k+1$ stand in the same column of the tableau $\La$.
\!Then the element $\Psi_\la\!\in\Sen$ \!is divisible on both sides by
$\phi_{k,\ts k+1}(z_k,z_{k+1})$.
\endproclaim

\nt
The proof of the next lemma is similar to that of Lemma 2.1 and will be
omitted.

\proclaim{Lemma 2.10}
Restriction of $\phik$ to the set of $(u,v,w)$
such that 
$$
\frac1{(v-w)^2}+\frac1{(v+w)^2}=1\,.
$$
is continuous at $u=\pm w$.
\endproclaim

\nt
We will conclude this section with two brief remarks. First of all,
notice that the proofs of Lemma 2.1 and Theorem 2.2 yield explicit formulas
for both elements $\Psi_\la,\Up_\la\in\Sen$. The element $\Up_\la\in\Sen$
should be regarded as an analogue of the classical Young symmetrizer
[\hs16\hs] in
the group ring $\CC\cdot S_n$; see [\hs11\hs] for more detail.

\kern15pt
\centerline{\bf3.\ Jucys-Murphy elements of the Sergeev algebra}\section{\,}
\kern-20pt

\nt
For each $k=1,2,\ldots$ we can regard the algebra $H_k$ as a subalgebra in
$H_{k+1}$.
Here the symmetric group $S_k\subset S_{k+1}$ acts on $1\lc k+1$
preserving the
number $k+1$. So we get a chain of subalgebras
$$
H_1\subset H_2\subset\ldots\subset H_n\subset H_{n+1}\subset\ldots\,\,.
$$
For each $k=1,2,\ldots$ introduce the element of the algebra $H_k$
$$
x_{k}=\sum_{1\leqslant i< k}\,(i\,k)+(i\,k)\,a_ia_{k}\,.
\Tag{2.86}
$$
In particular, $x_1=0$. The following proposition will be used in the next
section.

\proclaim{Proposition 3.1}
We have equality of rational functions of $u$ valued in $H_{n+1}$
$$
\prod^\rightarrow_{1\le i\le n}\ \phi_{n+1,i\,}(u\ts,z_i)
\cdot\Psi_\la=
\Bigl(\,1-\frac{x_{n+1}}u\,\Bigr)
\cdot\Psi_\la\ts.
\Tag{2.99}
$$
\endproclaim

\demo{Proof}
Denote by $\Xiu$ the rational function at the left hand side of \(2.99).
The value of this function at $u=\infty$~is~$\Psi_\la$. Moreover, the residue
of $\Xiu$ at $u=0$ is $-\,x_{n+1}\Psi_\la$.
It remains to prove that $\Xiu$
has a pole only at $u=0$ and this pole is simple. Let an index
$i\in\{2\lc n\}$ be fixed. The factor
$\phi_{n+1,i\,}(u,z_i)$
in \(2.99) has a pole at $u=\pm z_i$.
We shall prove that when estimating from
above the order of the pole at $u=\pm z_i$ of $\Xiu,$ that factor does not
count.

Suppose that the number $i$ does not appear in the first row of the tableau
$\La$.
Then the number $i-1$ appears in $\La$ straight above $i$.
In particular, $c_{i-1}=c_i+1$. By Corollary 2.8 the element $\Psi_\la$
is divisible on the left by
$$
\psi_{i-1,i}(c_{i-1}\ts,c_i)=\phi_{i-1,i}(z_{i-1}\ts,z_i).
$$
But the product
$$
\phi_{\ts n+1,i-1}(u\ts,z_{i-1})\ts\ts
\phi_{\ts n+1,i}(u\ts,z_{i})\ts\ts
\phi_{i-1,i}(z_{i-1}\ts,z_i)=
\phi_{\ts n+1,i-1,i\,}(u\ts,z_{i-1},z_i)
$$
is regular at $u=c_i$ due to Lemma 2.10. 

Now suppose that the number $i$ appears in the first row of $\La$.
Let $k$ be the number $k$ adjacent to $i$ on the left in the first row
of $\La$. Then $c_k=c_i-1$.
By Corollary~2.5 the element $\Psi_\la$ is divisible on the left by
$$
\prod_{k<j<i}^\leftarrow\ts
\phi_{\ts ji\ts}(z_j\ts,z_i)
\cdot
\phi_{\ts ki\ts}(z_k\ts,z_i).
\Tag{2.98}
$$
Consider the product of factors in \(2.99)
$$
\phi_{\ts n+1,k\ts}(u\ts,z_k)
\cdot\!\!
\prod_{k<j<i}^\rightarrow\!
\phi_{\ts n+1,j\ts}(u\ts,z_j)
\cdot
\phi_{\ts n+1,i\ts}(u\ts,z_i).
\Tag{2.97}
$$
Multiplying the product \(2.97) on the right by \(2.98) and using
\(2.1),\(2.2)
we get
$$
\prod_{k<j<i}^\leftarrow\ts
\phi_{\ts ji\ts}(z_j\ts,z_i)
\cdot
\phi_{\ts n+1,k,i\,}(u,z_k\,,z_i)\ts\ts
\cdot
\prod_{k<j<i}^\rightarrow\!
\phi_{\ts 1,j+1\ts}(u\ts,c_j).
$$
The latter product is regular at $u=z_i$ by Lemma 2.10. The proof is complete
\enddemos

\nt
The elements $x_1\lc x_n$ of the algebra $H_n$ will play an important role
in this
article. They will be called the {\it Jucys-Murphy elements\,} of the algebra
\text{$H_n$; cf.~[\hs17\hs].}

\proclaim{Lemma 3.2}
a) The elements $x_1\lc x_n$ of the algebra $H_n$ pairwise commute.
b) For any $r=1,2,\ldots$ the element $x_1^{2r}+\ldots+x_n^{2r}$ belongs to
the centre of $H_n$\,.
\endproclaim

\demo{Proof}
The first statement of this lemma can be verified by direct calculation. One
can also check directly the relations in the algebra $H_n$
$$
\gather
a_k\,x_k=-\,x_k\,a_k; 
\qquad
a_k\,x_l=x_l\,a_k\,,
\quad
k\neq l\,;
\Tag{2.121212}
\\
(k,k+1)\,x_{k+1}-x_k\,(k,k+1)=1+a_ka_{k+1}\,;
\Tag{2.121213}
\\
x_{k+1}\,(k,k+1)-(k,k+1)\,x_k=1-a_ka_{k+1}\,.
\Tag{2.121214}
\endgather
$$
By making use of \(2.121213) and \(2.121214) one can verify
$$
\bigl[\,(k,k+1)\hskip1pt,x_k\,x_{k+1}\,\bigr]=0\,,
\quad
\bigl[\,(k,k+1)\hskip1pt,x_k^2+x_{k+1}^2\,\bigr]=0\,
\Tag{2.121215}
$$
where the square brackets stand for the commutator. By the definition \(2.86)
the transposition $(k,k+1)\in S_n$ also commutes with $x_l$
if $l\neq k,k+1$. So the commutation relations \(2.121212),\(2.121215) imply
the second statement of Lemma 3.2
\enddemos

\nt
According to the following theorem, the element $\Psi_\la\in H_n$ is a joint
eigenvector of $x_1\lc x_n$
with respect to the left regular action of the algebra $H_n$.

\proclaim{Proposition 3.3}
For each $k=1\lc n$ we have $x_k\cdot\Psi_\la=z_i\cdot\Psi_\la$.
\endproclaim

\demo{Proof}
Note that $x_1=0$ and $z_1=0$ by definition. So we will assume that
$k\neq1$. For any $k=2\lc n$
an argument similar to that already used in the proof of Proposition 3.1 
shows the equality of rational functions of $u$ valued in $H_n$
$$
u\,\cdot\!\!
\prod^\rightarrow_{1\le i< k}\ \phi_{ki}(u\ts,z_i)
\cdot\Psi_\la=
(u-x_k)
\cdot\Psi_\la\ts.
\Tag{2.999}
$$
Denote by $\Phi_k(u)$ the rational function appearing at the
left hand side of \(2.999). We shall prove that $\Phi_k(z_k)=0$.

Consider any factor $\phi_{ki}(u\ts,z_i)$
at the left hand side of \(2.999) with $z_i=z_k$.
Then this factor has a pole at $u=z_k$. However,
according to the proof of Proposition 3.1 this factor
does not count when estimating from
above the order of the pole at $u=z_k$ of $\Phi_k(u),$ provided
$i\neq1$.
But $u\cdot\phi_{k1}(u\ts,z_1)$ has no pole at $u=0$.

Suppose the number $k$ does not appear in the first row of the tableau $\La$.
Then the number $k-1$ appears in $\La$ straight above $k$. So
$c_{k-1}=c_k+1$.
Consider the factor
$\phi_{k,k-1}(u\ts,z_{k-1})$ in \(2.999).
By Corollary 2.8 the element $\Psi_\la$
is divisible on the left by
$$
\psi_{k-1,k}(c_{k-1}\ts,c_k)=\phi_{k-1,k}(z_{k-1}\ts,z_k).
$$
But the product
$$
\phi_{k,k-1}(u\ts,z_{k-1})\cdot\phi_{k-1,k}(z_{k-1}\ts,z_k)
$$
vanishes at $u=z_k$ due to \(2.31). So the function 
$\Phi_k(u)$ vanishes at $u=z_k$ as well.

Now suppose that the number $k$ appears in the first row of $\La$. Note that
$k\neq 1$ by our assumption.
Let $i$ be the number adjacent to $k$ on the left in the first row of $\La$.
Then $c_i=c_k-1$. 
Consider the product of the factors in \(2.999)
$$
\phi_{\ts ki\ts}(u\ts,z_i)
\cdot\!\!
\prod_{i<j<k}^\rightarrow\!
\phi_{\ts kj\ts}(u\ts,z_j)
\Tag{2.998}
$$
None of them has a pole at $u=z_k$.
By Corollary~2.5 the element $\Psi_\la$ is divisible on the left by
$$
\prod_{i<j<k}^\leftarrow\ts
\phi_{\ts jk \ts}(z_j\ts,z_k)
\cdot
\phi_{\ts ik\ts}(z_i\ts,z_k).
$$
Due to \(2.31) when we multiply the latter product on the left by \(2.998)
and set $u=z_k$~we~get 
$$
\prod_{i\le j<k}\ts
\biggl(1-\frac1{(z_j-z_k)^2}-\frac1{(z_j+z_k)^2}\biggr)
$$
which is equal to zero because the factor corresponding to $j=i$ vanishes
\enddemos

\nt
Let $\mu$ run through the set $\Cal D_\la$
of all strict partitions of $n-1$ which can
be obtained by decreasing one of the parts of $\la$. Put $m_{\la\mu}=1$
if the number $\ell_\mu$ is odd while $\ell_\la$ is even; otherwise put
$m_{\la\mu}=2$. Restriction  of
the $\Sen$-module $U_\la$ to $H_{n-1}$ splits into a direct sum of
the modules
$U_\mu$ where each $U_\mu$ appears $m_{\la\mu}$ times;
see [\hs9\hs,\hs Lemma 6\hs].
This branching property determines the irreducible module $U_\la$
over $\ZZ_2$-graded algebra $\Sen$ uniquely. 
We set $H_0=\CC$.

\proclaim{Theorem 3.4}
Under the left regular action of $H_n$
\text{the space $\Sen\cdot\Psi_\la$} splits into a direct sum of copies of
the irreducible module $U_\la$.
\endproclaim

\demo{Proof}
Let us prove by induction on $n$ the following statement: for any
$r=1,2,\ldots$
the central element $x_1^{2r}+\ldots+x_n^{2r}\in\Sen$ acts as
$$
z_1^{\hskip1pt 2r}+\ldots+z_n^{\hskip1pt 2r}\,=\,
c_1^{\,r}(c_1+1)^r+\ldots+c_n^{\,r}(c_n+1)^r
\Tag{2.999999999}
$$
in the irreducible module $U_\la$. Due to Proposition 3.3 the latter statement
implies Theorem~3.4. Indeed, the collection of the
eigenvalues \(2.999999999) determines the shifted Young diagram
$\la$ uniquely. Therefore if an irreducible $\Sen$-module $U$  
contains an eigenvector of $x_1^2\lc x_n^2$ with the respective eigenvalues
$z_1^2\lc z_n^2$ then $U=U_\la$\,.

If $n=1$ then $x_1=0$ and $c_1=0$ so the statement to prove is trivial.

Now assume that $n>1$. Suppose there exist two distinct diagrams
$\mu,\nu\in\Cal D_\la$. Each of them can be obtained by adding
a box to the same shifted diagram.
Let $c$ and $d$ be the contents of these two boxes.
Consider the irreducible components $U_\mu$ and $U_{\nu}$ in the restriction
of $U_\la$ to $H_{n-1}$. The element $x_n\in\Sen$ commutes with the 
subalgebra $H_{n-1}$ and acts in $U_\mu$,$U_{\nu}$ by certain
numbers $z,w\in\CC$ respectively. By comparing the actions
of the central elements
$x_1^{2r}+\ldots+x_n^{2r}\in\Sen$ in $U_\mu$,$U_{\nu}$ and 
by applying the inductive assumption to these $H_{n-1}$-modules
we get the equality
$$
c^{\,r}(c+1)^r+z^{2r}=d^{\,r}(d+1)^r+w^{2r}\,.
$$
These equalities for $r=1,2$ imply that $z^2=d(d+1)$ and $w^2=c(c+1)$.
So we obtain the required statement.

Assume that $\Cal D_\la=\{\hskip1pt\mu\hskip1pt\}$. 
Take a joint eigenvector $\xi\in U_\la$ of $x_n,x_{n-1}\in H_n$.
Let $u,v\in\CC$ be the respective eigenvalues.
Observe that $u\neq\pm\,v$. If $u=v$ then
by applying \(2.121213),\(2.121214) with $k=n-1$ to the vector $\xi$ we
obtain that
$(x_n-x_{n-1})\cdot(k,k+1)\,\xi=2\,\xi$ while $(x_n-x_{n-1})\cdot\xi=0$.
This contradicts the property $(x_n-x_{n-1})^\ast=x_n-x_{n-1}$.
Thus $u\neq v$. Similarly,
by taking the vector $a_{n-1}\cdot\xi$ instead of $\xi$ we prove that
$u\neq-\,v$. Now consider the element
$$
h_n=(n-1,n)\cdot(x_{n-1}^2-x_n^2)+(x_{n-1}+x_n)-
a_{n-1}\,a_n\,(x_{n-1}-x_n)\in\Sen\,.
$$
One can derive that $x_n\,h_n=h_n\,x_{n-1}\,$ from \(2.121214). So 
$x_n\cdot h_n\,\xi=v\cdot h_n\,\xi$. But $x_n\,c_n=-c_n\,x_n$ while
$d_{\la\mu}\le2$.
So any eigenvalue of $x_n$ in $U_\la$ is either $u$ or $-\ts u$.
Hence $h_n\cdot\xi=0$, that is
$$
(n-1,n)\cdot\xi=\Bigl(\,\frac1{u-v}+\frac{a_{n-1}\,a_n}{u+v}\,\Bigr)
\cdot\xi\,.
\Tag{2.001}
$$
In particular, the pair $(u,v)$ satisfies the condition \(2.34)
since $(n-1,n)^2\cdot\xi=\xi$. 

Let us now take any $H_{n-1}$-irreducible component $U_\mu\subset U_\la$.
We can choose $\xi\in U_\mu$. We will show that $u^2=z_n^2$. Then we will
get the required statement by the inductive assumption.
If $\la=(2)$ then $v=0$, so $u^2=2=z_2^2$
by \(2.34). We will assume that $n>2$. Let us choose $\xi\in U_\mu$ so that
$x_i\cdot\xi=z_i\cdot\xi$ for $i=n-1,n-2$. Here we use Proposition 3.3
and the inductive assumption. So we have $v=z_{n-1}$. We will also write
$w=z_{n-2}$. Consider the following five cases. 

\smallskip\nt
i) Suppose that $\ell_\la=2$ while both parts of $\la$ are greater than 1.
Then the number $n-1$ stands straight above $n$ in $\La$ while
$n-2$ stands next to the left of $n$. Now
$$
\gather
x_{n-1}\cdot h_{n-1}\,\xi=
h_{n-1}\,x_{n-2}\cdot\xi=
w\cdot h_{n-1}\,\xi\,,
\\
h_{n-1}\,\xi=
(w^2-v^2)\,\Bigl(
(n-2,n-1)-\frac1{v-w}-
\frac{a_{n-2}\,a_{n-1}}{v+w}
\,\Bigr)\cdot\xi\neq0.
\endgather
$$
So the pair $(u,w)$ satisfies the condition \(2.34) as well the pair $(u,v)$. 
Note that here $w^2=(c_n-1)\,c_n$ and $v^2=(c_n+1)(c_n+2)$.
Therefore $u^2=c_n\,(c_n+1)=z_n^2$.

In the remaining four cases the number $n-1$
will stand next to $n-2$ in the tableau $\La$. Due to Corollaries 2.5 and 2.9
we can now assume that
$$
(n-2,n-1)\cdot\xi=\Bigl(\,\frac1{v-w}+\frac{a_{n-2}\,a_{n-1}}{v+w}\,\Bigr)
\cdot\xi\,.
\Tag{2.002}
$$

\smallskip\nt
ii) Suppose that $\ell_\la>2$. Then the numbers $n-1,n-2$  stand straight
above $n$
in the tableau $\La$. So $v^2=(c_n+1)(c_n+2)$ and $w^2=(c_n+2)(c_n+3)$.
In particular, here we have $w\neq0$. The condition
\(2.34) implies that either $u^2=c_n\,(c_n+1)$ or $u^2=(c_n+2)(c_n+3)=w^2$.
But if $u=\pm\,w$ then by substituting \(2.001),\(2.002) in
$$
(n-1,n)\,(n-2,n-1)\,(n-1,n)\cdot\xi=(n-2,n-1)\,(n-1,n)\,(n-2,n-1)\cdot\xi
$$
we obtain that $w=0$. This contradiction demonstrates that
$u^2=c_n\,(c_n+1)=z_n^2$.

\smallskip\nt
iii) If $\la=(2,\!1)$ then $m_{\la\mu}=1$, so $u=-\,u$.
Thus $u^2=0=z_3^2$. Note that here
$$
(12)\cdot\xi=\frac1{\sqrt2}\,(a_1a_2+1)\cdot\xi\,,
\qquad
(23)\cdot\xi=\frac1{\sqrt2}\,(a_2\,a_3-1)\cdot\xi\,.
\Tag{2.003}
$$

\smallskip\nt
iv) Suppose that $\ell_\la=1$ and $n>3$. Then  $n-1,n-2$ stand straight
to the left of $n$ in the tableau $\La$. So
$v^2=(c_n-1)\,c_n$ and $w^2=(c_n-2)(c_n-1)$. Now
\(2.34) implies that either $u^2=c_n\,(c_n+1)$ or $u^2=(c_n-2)(c_n-1)=w^2$.
But $w\neq0$. As in (ii) one proves that only the
case $u^2=c_n\,(c_n+1)=z_n^2$ is possible.
\vbox{
\smallskip\nt
v) Finally, suppose that $\la=(3)$. Here $v^2=2$ and \(2.34)
implies that either $u^2=6$ or $u^2=0$. But in the latter case the
transpositions
$(12),\!(23)\in H_3$ would act in the irreducible module $U_{(3)}$
by the same formulas \(2.003) as in the module $U_{(2,1)}$.
Since the $H_3$-modules $U_{(3)}$ and $U_{(2,1)}$ are non-equivalent,
we get $u^2=6=z_3^2\quad\square$
}
\enddemo

\nt
More detailed analysis shows that the left ideal $\Sen\cdot\Psi_\la$
splits into a direct sum of $2^{[\ell_\la/2]}$ copies of the irreducible 
\text{module $U_\la$ [\hs11\hs,\hs Theorem 8.3\hs].} However,
we will not use the latter fact in the present article. We will need the
following corollary to Theorem 3.4. Let $h$ run through the set of
basic elements of $\Sen$
$$
g\,a_1^{l_1}\ldots a_n^{l_n}\,;
\qquad
g\in S_n\,,
\quad
l_1\lc l_n=0,1.
\Tag{5.3125}
$$
We denote by $\chi_\la(h)$ the trace of the element $h\in\Sen$ in the
module $U_\la$ normalized so that $\chi_\la(1)=1$.
Consider the corresponding central element of $\Sen$
$$
\Chi_\la=\sum_h\ \chi_\la(h)\cdot h\1\,.
\Tag{5.625}
$$

\proclaim{Corollary 3.5}
We have the equality
$$
\Chi_\la\cdot{2^n\,n\,!\,}=
\sum_h\ h\,\,\Psi_\la\,h\1\,.
$$
\endproclaim

\nt
Now fix any non-negative integer $m$ and consider $\Sen$
as a subalgebra in $H_{n+m}$. 
We will write $\cj=j+n$ for every $j=1\lc m$.
For each $k=1\lc n$ denote
$$
y_{k}=\!\!\sum_{1\leqslant j\leqslant m}\,
(\,k\,\cj\,)-(\,k\,\cj\,)\,a_{k}\,a_{\cj}\,.
\Tag{5.1}
$$
The next auxiliary statement easily follows from Proposition 3.3;
\text{cf.~[\hs13\hs].}

\proclaim{Proposition 3.6}
The product $\Psi_\la\cdot(y_1-z_1)\ldots(y_n-z_n)\in H_{n+m}$ equals
the sum
$$
\sum_{\ j_1\ldots j_n}\,
(\,1\,\cj_1\,)\ldots(\,n\,\cj_n\,)
\cdot
(1-a_1a_{\cj_1})\ldots(1-a_na_{\cj_n})
\Tag{5.3}
$$
where all the indices $j_1\lc j_n\in\{1\lc m\}$ are pairwise distinct.
\endproclaim

\demo{Proof}
We will use the induction on $n$. Note that $z_1=0$ so for $n=1$
we get the required statement by the definition \(5.1). Let us now
suppose that $n>1$.
Then by the inductive assumption and by the definition \(5.1)
we have the equality
$$
\align
\Psi_\la\cdot(y_1-z_1)\ldots(y_n-z_n)&=
\Tag{5.2}
\\
\Psi_\la\,\cdot\!\!\!\!\sum_{j_1\ldots j_{n-1}}\!\!
(\,1\,\cj_1\,)\ldots(\,n-1,\cj_{n-1}\,)
\,
\bigl(1-a_1a_{\cj_1}\bigr)&\ldots\bigl(1-a_{n-1}a_{\cj_{n-1}}\bigr)
\,\x
\\
\x\,\sum_j\,(\,n\,\cj\,)\,\bigl(1-a_{n}\,a_{\cj}\,\bigr)&\ -
\\
\Psi_\la\,\cdot\!\!\!\!\sum_{j_1\ldots j_{n-1}}\!\!
z_n\,(\,1\,\cj_1\,)\ldots(\,n-1,\cj_{n-1}\,)
\,
\bigl(1-a_1a_{\cj_1}\bigr)&\ldots\bigl(1-a_{n-1}a_{\cj_{n-1}}\bigr)
\endalign
$$
where $j_1\lc j_{n-1},j\in\{1\lc m\}$ and the indices $j_1\lc j_{n-1}$
are pairwaise distinct.
Due to Lemma 2.3 and Proposition 3.3 by the definition \(2.86)
the sum in the last line of this equality can be replaced by the sum
$$
\sum_{i}\sum_{j_1\ldots j_{n-1}}
(n\,i)\,\bigl(1-a_na_i\bigr)\!\cdot
(\,1\,\cj_1\,)\ldots(\,n-1,\cj_{n-1}\,)
\,
\bigl(1-a_1a_{\cj_1}\bigr)\ldots\bigl(1-a_{n-1}a_{\cj_{n-1}}\bigr)
$$
where the index $i$ runs through $1\lc n-1$. The latter sum can be rewritten
as
$$
\sum_{j_1\ldots j_{n-1}}\sum_{j}\ \ 
(\,1\,\cj_1\,)\ldots(\,n-1,\cj_{n-1}\,)
\cdot(\,n\,\cj\,)\,\bigl(1-a_na_{\cj}\bigr)\cdot
\bigl(1-a_1a_{\cj_1}\bigr)\ldots\bigl(1-a_{n-1}a_{\cj_{n-1}}\bigr)
$$
where $j$ runs through $j_1\lc j_{n-1}$. Now Proposition 3.6
follows from \(5.2)
\enddemos

\nt
Let $\mu$ be any strict partition of $m$. We will 
identify the partitions $\la$ and $\mu$ with their
shifted Young diagrams. Take the embedding of the algebra $H_m$
into $H_{n+m}$ such that the symmetric group $S_m\subset S_{n+m}$ acts by
permutations of the numbers $n+1\lc n+m$. Denote by
$\Psi_\mu\!\!\!\tilde{\phantom{o}}$ the image of the element
$\Psi_\mu\in H_m$
with respect to this embedding.

\proclaim{Corollary 3.7}
We have the equality
$\Psi_\la\cdot(y_1-z_1)\ldots(y_n-z_n)\cdot
\Psi_\mu\!\!\!\tilde{\phantom{o}}=0$
if the diagram $\la$ is not contained~in~$\mu$.
\endproclaim

\demo{Proof}
If $m<n$ there is no summand in \(5.3) and
$\Psi_\la\cdot(y_1-z_1)\ldots(y_n-z_n)=0$.
Suppose that $m\ge n$ but the diagram $\la$ is not contained in $\mu$.
Consider $H_n$ as a subalgebra in  $H_m$ with respect to the standard
embedding.
The restriction of the $H_m$-module $U_\mu$ to $H_n$ does not contain
any irreducible component isomorphic to $U_\la$. So by Lemma 2.3 and
Theorem 3.4 we get $\Psi_\la\cdot\Psi_\mu=0$ in $H_m$. Moreover,
$\Psi_\la\cdot a_{j_1}\ldots a_{j_s}\!\cdot\Psi_\mu=0$
for any $j_1\lc j_{s}\in\{1\lc m\}$. Now the required equality follows
from Proposition 3.6
\enddemos

\nt
In the next section we interpret Corollary 3.7 in terms of 
classical invariant theory.


\kern15pt
\centerline{\bf 4. Capelli identity for the queer Lie superalgebra}
\section{\,}
\kern-20pt

\nt
In this section we will let the indices $i,j$ run through $\pm\,1\lc\pm\,N$.
We will write $\bi=0$ if $i>0$ and $\bi=1$ if $i<0$. Consider the
$\ZZ_2$-graded vector space $\CN$. Let $e_i\in\CN$ be an element of the
standard basis. The $\ZZ_2$-gradation on $\CN$ is defined so that
$\deg e_i=\bi$.
Let $E_{ij}\in\EndCN$ be the standard matrix units. The algebra $\EndCN$
is $\ZZ_2$-graded so that $\deg E_{ij}=\bi+\bj$.

We will also regard $E_{ij}$ as generators of the complex Lie superalgebra
$\glN$.
The {\it queer} Lie superalgebra $\qN$ is a subalgebra in $\glN$ spanned by
the
elements $F_{ij}=E_{ij}+E_{-i,-j}$. Thus $F_{-i,-j}=F_{ij}$ by definition.
Note that the image of the defining representation $\qN\to\EndCN$ coincides
with
the supercommutant of the element
$$
J=\sum_{j}\ts E_{j,-j}\cdot{(-1)}^{\ts\bj}\in\EndCN.
$$ 

In this section we will use the following convention.
Let $\operatorname{A}$ and $\operatorname{B}$ 
be any~two associative complex
$\ZZ_2$-graded algebras. Their tensor product
$\operatorname{A}\ot\operatorname{B}$ will be
a $\ZZ_2$-graded algebra such that for any homogeneous elements
$X,X^\prime\in\operatorname{A}$ and $Y,Y^\prime\in\operatorname{B}$
$$
\align
(X\ot Y)\ts (X^\prime\ot Y^\prime)&=X\ts X^\prime\ot Y\ts Y^\prime
\cdot(-1)^{\ts\deg X^\prime\deg Y},
\\
\deg\ts(X\ot Y)&=\deg X+\deg Y\ts.
\endalign
$$
If the algebra $\operatorname{A}$ is unital denote
by $\iota_s$ its embedding into the tensor product
$\operatorname{A}^{\!\ot n}$
as the $s$-\ts th tensor factor:
$$
\iota_s(X)=1^{\ot\ts (s-1)}\ot X\ot1^{\ot\ts(n-s)}\,,
\qquad
1\le s\le n.
$$
We will also use various embeddings of the algebra
$\operatorname{A}^{\!\ot\ts m}$ into $\operatorname{A}^{\!\ot\ts n}$
for any $m\leqslant n$.
For any choice of pairwise distinct indices $s_1\lc s_m\in\{\ts1\lc n\ts\}$
and an element $X\in\operatorname{A}^{\!\ot m}$ of the form
$X=X^{(1)}\ot\ldots\ot X^{(m)}$ we will denote
$$
X_{s_1\ldots s_m}=\ts
\iota_{s_1}\bigl(X^{(1)}\bigr)\ts\ldots\,\ts\iota_{s_m}\bigl(X^{(m)}\bigr)
\in\operatorname{A}^{\!\ot n}.
$$

Denote 
$$
P=\sum_{i,j}\ts E_{ij}\ot E_{ji}\cdot{(-1)}^{\ts\bj}\in\EndCN^{\,\ot2}
$$
For any positive integer $n$ a representation $\Sen\to\EndCN^{\ot\ts n}$
of the Sergeev algebra can be determined by the assignments
$$
a_k\mapsto J_k
\quad\text{and}\quad
(\ts k\ts l\ts)\mapsto P_{kl}\ts.
\Tag{3.55555}
$$

Let us now consider the enveloping algebra $\UN$ of the Lie superalgebra
$\qN$.
The algebra $\UN$ is a Hopf superalgebra: the comultiplication, counit and
the
antipodal map are defined for $F\in\qN$ respectively by
$$
\Delta\ts(F)=F\ot1+1\ot F\ts,
\ \quad
\ep(F)=0\ts,
\ \quad
S(F)=-\,F\ts.
$$

Take the $\!n$-th tensor power 
of the defining representation $\UN\to\EndCN$.
Its image coincides [\hs9\hs,\hs Theorem 3\hs] with supercommutant
of the image of Sergeev algebra $\Sen$ relative to \(3.55555).
Let $\la$ be any strict partition of $n$.
Take the irreducible module $U_\la$ over the $\ZZ_2$-graded
algebra $\Sen$. Denote by $V_\la$ the subspace in
$\operatorname{Hom}\hskip1pt\bigl(U_\la\,,\EndCN^{\ot\ts n}\bigr)$
consisting of all the elements which supercommute with the action of $\Sen$.
We have an irreducible representation
of the $\ZZ_2$-graded algebra $\UN$ in the space $V_\la$.
We denote it by $\pi_\la$.
Here $V_\la\neq\{\ts0\ts\}$ if and only if $\ell_\la\le N$
[\hs 9\hs,\hs Theorem 4\hs].
From now on we will assume that this is the case.

Now let two positive integers $N$ and $M$ be fixed.
Let the indices $a,b$ run through $\pm\,1\lc\pm\,M$ while
the indices $i,j$ keep running through $\pm\,1\lc\pm\,N$
We will also write $\ba=0$ if $a>0$ and $\ba=1$ if $a<0$.
We will use the generators $F_{ij}\in\qN$ and $F_{ab}\in\qM$.

Introduce a supercommutative algebra $\PMN$ with the free generators
$x_{ia}$ where $a>0$ and the $\ZZ_2$-gradation is defined by 
$\deg x_{ia}=\bi$.
It will be convenient to set
$$
x_{i,-a}=\sqrt{-1}\cdot x_{-i,a}\,;
\quad
a>0.
$$
Let $\d_{ia}$ be the left derivations in the supercommutative algebra $\PMN$
corresponding to the generators $x_{ia}$. Here we allow both indices $i$ and
$a$
to be negative, so that
$$
\d_{i,-a}=-\,\sqrt{-1}\cdot\d_{-i,a}\,;
\quad
a>0.
$$
The algebra generated by all the left derivatives $\d_{ia}$ in $\PMN$ along
with
the operators of left multiplication by $x_{jb}$ will be denoted by $\PDMN$.
Note that for the arbitrary
indices $i=\pm\,1\lc\pm N\,$ and $a=\pm\,1\lc\pm\,M$
we have $\deg x_{ia}=\bi+\ba$ in $\PMN$.

\proclaim{Proposition 4.1}
The assignments for $i,j=\pm1\lc\pm N$ and $a,b=\pm1\lc\pm M$
$$
F_{ij}\mapsto\sum_b\,x_{ib}\,\d_{jb}
\quad\ \text{and}\ \quad
F_{ab}\mapsto\sum_j\,x_{ja}\ts\d_{jb}\cdot(-1)^{\,\bj\,(\ba\,+\,\bb)}
\Tag{7.1}
$$
define representations in $\PDMN$ of the Lie superalgebras $\qN$ and $\qM$.
The images of 
$\UN$ and $\UM$ in these representations are the supercommutants of
each~other. 
\endproclaim

\demo{Proof}
Let $e_i\in\CN$ and $e_a\in\CM$ be elements
of the standard bases. Identify 
the tensor product $\End(\CN)\ot\End(\CM)$
with the algebra $\End(\CMN)$ so that
$$
E_{ij}\ot E_{ab}\cdot e_k\ot e_c=e_i\ot e_a
\cdot\de_{jk}\,\de_{bc}\,(-1)^{\,\bj\,(\ba\,+\,\bb)}\,.
$$
Now the standard embeddings 
$X\mapsto X\ot1$ and $Y\mapsto1\ot Y$ of
$\End(\CN)$ and $\End(\CM)$ into $\End(\CN)\ot\End(\CM)$ respectively
define actions of the Lie superalgebras $\qN$ and $\qM$ in the space $\CMN$. 
These actions preserve the subspace $W$ in $\CMN$ with the basis
consisting of the vectors
$$
e_i\ot e_a-\sqrt{-1}\cdot e_{-i}\ot e_{-a}\,;
\quad
a>0.
$$
The supersymmetric algebra $\operatorname{S}(\hskip1pt W)$ coincides with
$\PMN$,
the above basic vector~of~$W$ being identified with the generator $x_{ia}\,$.
The action of $\qN$ and $\qM$ in the $\operatorname{S}(\hskip1pt W)$ is then
given by \(7.1). So the images of the enveloping algebras
$\UN$ and $\UM$ in $\PDMN$ supercommute. Moreover, due to
[\hs9\hs,\hs Theorem 3\hs] the spectrum of the $\UN\ot\UM$-module
$\operatorname{S}(\hskip1pt W)$ is simple.  Hence the images
of $\UN$ and $\UM$ in $\PDMN$ are the supercommutants of each other
\enddemos

\nt
Thus we have defined representations $\UN\to\PDMN$ and $\UM\to\PDMN$,
we will denote them by $\ga$ and $\ga^{\,\prime}$ respectively.
In this section we will consider the subspace
$$
\I=
\ga\hskip1pt\bigl(\hskip1pt\UN\bigr)\cap
\ga^{\,\prime}\bigl(\hskip1pt\UM\bigr)\subset\PDMN.
$$
Next corollary to Proposition 4.1 gives another description of
the subspace~$\I\subset\PDMN$.

\proclaim{Corollary 4.2}
The subspace $\I\subset\PDMN$ consists of those operators
which {supercom\-mute} with the images of $\ga$ and $\ga^{\,\prime}$.
\endproclaim

\nt
Let us also equip the algebra $\PMN$ with the $\ZZ$-gradation such that
$\deg x_{ia}=1$. Denote by $\PMN_n$ the subspace in $\PMN$ consisting
of the elements of degree $n$. According to [\hs9\hs,\hs Theorem 3\hs] 
the irreducible components of the $\UN\ot\UM$-module
$\PMN_n$ 
are parametrized by the strict
partitions $\la$ of $n$ with not more than $M,N$ parts.
Denote by $W_\la$ the irreducible component of $\PMN_n$ corresponding to
$\la$. If the number $\ell_\la$ is even then
$W_\la=V_\la\ot V_\la^{\hskip1pt\prime}$
where the second tensor factor is the irreducible $\UM$-module
corresponding to the partition $\la$. If $\ell_\la$ is odd then the tensor
product $V_\la\ot V_\la^{\hskip1pt\prime}$ splits into direct sum of
two copies of the irreducible module $W_\la$.

Let $\mu$ run through the set of
strict partitions of $m=0,1,2,\ldots$ with not more than $M,N$ parts.
The next proposition gives a distinguished decomposition of the vector space
$\I$ into a direct sum of one-dimensional subspaces;
\text{cf. [\hs3\hs,\hs Theorem 1\hs].}

\proclaim{Proposition 4.3}
There exists a unique one-dimensional subspace $\I_\la\in\I$
such that $\I_\la\cdot W_\mu=\{0\}$ when $m<n$ or $m=n$ but $\mu\neq\la$.
We have a decomposition
$$
\I=\underset\la\to\oplus\ \I_\la\,
\Tag{7.0}
$$
where $\la$ runs through the set of all strict partitions with not more than
$M,N$ parts.
\endproclaim

\demo{Proof}
Consider the actions $\operatorname{ad}\!\crc\!\ga$ and
$\operatorname{ad}\!\crc\!\ga^{\,\prime}$ in the space $\PDMN$ of the
Lie superalgebras $\qN$ and $\qM$ respectively. By Corollary 4.2
the subspace $\I\subset\PDMN$ consists of all the invariants of
both actions. Denote by $\DMN$ the subalgebra in $\PDMN$ generated
by all the left derivatives $\d_{ia}$. The actions of $\qN$ and $\qM$
in $\PDMN$ preserve the subspace $\DMN$. Obviously, we have
$\PDMN=\PMN\cdot\DMN$. Moreover, 
as a $\UN\ot\UM$-module the space $\PDMN$ now decomposes into the tensor
product $\PMN\ot\DMN$.

Let us equip the algebra $\DMN$ with the $\ZZ$-gradation such that
$\deg\d_{ia}=1$. Denote by $\DMN_n$ the subspace in $\DMN$ consisting
of the elements of degree $n$. As a module over $\UN\ot\UM$ the subspace
$\DMN_n$ splits into direct sum of irreducible modules $W_\la^\ast$
such that the tensor product $W_\la\ot W_\mu^\ast$ contains an
invariant subspace if and only if $\la=\mu$.
Let $\I_\la$ be the invariant subspace in $W_\la\ot W_\la^\ast\subset\PDMN$,
it is one-dimensional.

Choose any non-zero element $I_\la\in\I_\la$.
Since $W_\mu\subset\PMN_m$ we have $I_\la\cdot W_\mu=0$ for $m<n$.
Suppose that $m=n$. Consider the linear map
$\PMN_n\to\PMN_n:\,X\to I_\la\cdot X$. It commutes with the actions of
$\qN$ and $\qM$ in $\PMN_n$. But the image of this map is contained in
$W_\la$
by the definition of the space $\I_\la$. So the restriction of this map to
$W_\mu\subset\PMN_n$ with $\mu\neq\la$ is zero. 
Thus we obtain the decomposition \(7.0). It is unique
since the spectrum of the $\UN\ot\UM$-module $\PMN$ is simple
\enddemos 

\nt
Let $\ZN$ be the centre of the enveloping algebra $\UN$. By definition,
an element of $\UN$ is central if it supercommutes with any element
of $\UN$. However, the centre $\ZN$ consists of even elements only
[\hs18\hs,\hs Theorem 1\hs]. 
Note that the representation $\ga:\UN\to\PDMN$ is faithful when $M\ge N$.
Then we get the equality $\I=\ga\,\bigl(\ZN\bigr)$ by
Proposition 7.1\,; cf. [\hs19\hs,\hs Section 3\hs].
In this section for any $M,N$ we will give an explicit formula for a
non-zero element in $\I_\la$.
We will also construct a non-zero element in
$\ga\1(\I_\la)\cap\ZN$. An element of $(\ga^{\,\prime})\1(\I_\la)\cap\ZM$
can be constructed in a similar way.
Thus for any $M,N$ we will get an evidential proof of the
equality
$$
\I=\ga\,\bigl(\hskip1pt\ZN\bigr)=\ga^{\,\prime}\,\bigl(\hskip1pt\ZM\bigr).
$$

Now let $\la$ be any strict partition of $n$ with $\ell_\la\le N$.
Let $R_\la\in\EndCN^{\ot\ts n}$ correspond
to $\Psi_\la\in\Sen$ with respect to \(3.55555).
Let $L_\la\subset(\CN)^{\ot\ts n}$ be the image of
the endomorphism $R_\la$.
The representation of $\UN$ in the space $L_\la$
is a direct sum of $2^{\,[\ell_\la/2]}$
copies of the irreducible representation in $V_\la$;
see Theorem 3.4. Denote by $\om_\la$ the respective homomorphism
$\UN\to\End(L_\la)$.
We will identify the algebra $\End(L_\la)$ with
the subalgebra in $\EndCN^{\ot n}$ consisting of all the elements
which have the form $X\,R_\la=R_\la\,Y$ for some $X,Y\in\EndCN^{\ot n}$.

Denote
$$
F=\sum_{i,j}\,E_{ij}\ot F_{ji}\cdot(-1)^{\,\bj}\in\EndCN\ot\UN.
$$
Note that for any element $X\in\qN$ we have the equality
$$
\bigl[\,X\ot1+1\ot X\,,\,F\,\bigr]=0
\Tag{7.2}
$$ 
where the square brackets stand for the supercommutator. For $s=1\lc n$
we will write
$$
F_s=\iota_s\ot\id\,\,(F)\in\End(\CN)^{\,\ot n}\ot\UN.
$$
As well as in the previous section put $z_s=\sqrt{c_s(c_s+1)}$ where
$c_s$ is the content of the box with number $s$ in the
shifted column tableau $\La$. Consider the element
$$
F_\la=R_\la\ot1\cdot(F_1-z_1)\ldots(F_n-z_n)\in\End(\CN)^{\,\ot n}\ot\UN.
\Tag{7.55555}
$$

Now let $\mu$ by any strict partition of $m$ with not more than $N$ parts.
The next proposition makes the central part of this section\,;
cf. [\hs20\hs].

\proclaim{Proposition 4.4}
\!We have $\id\ot\pi_\mu(F_\la)=0$ if the diagram $\la$
is not contained~in~$\mu$. 
\endproclaim 

\demo{Proof}
We can replace in Proposition 4.4 the representation $\pi_\mu$ by
the direct sum $\om_\mu$ of its copies. With respect
to the defining representation $\UN\to\EndCN$
$$
\EndCN\ot\UN\to\EndCN^{\,\ot2}:\,F\mapsto P\,\bigl(1-J_1J_2\bigr). 
$$
So the element
$\id\ot\om_\mu(F_\la)$ of $\EndCN^{\,\ot n}\ot\End(L_\mu)
\subset\EndCN^{\,\ot\,(n+m)}$ equals
$$
R_\la\ot1\cdot
\prod_{1\le s\le n}^\rightarrow
\Bigl(\sum_{1\le r\le m}\!\!
P_{s,n+r}\bigl(1-J_s\,J_{n+r}\bigr)-z_s\Bigr)\ \cdot1\ot R_\mu\,.
$$
This product is the image in $\EndCN^{\,\ot\,(n+m)}$
of the element from the Sergeev algebra $H_{n+m}$
$$
\Psi_\la\cdot(y_1-z_1)\ldots(y_n-z_n)\cdot
\Psi_\mu\!\!\!\tilde{\phantom{o}}\,\,;
$$
see the end of the previous section.
By Corollary 3.7 the latter product vanishes if the shifted diagram
$\la$ is not contained in the shifted diagram $\mu$
\enddemos

\nt
Consider the linear functional $\tr:\EndCN^{\,\ot n}\to\CC$ called the
{\it supertrace}. By definition, we have
$$
\tr:\,E_{i_1j_1}\ot\ldots\ot E_{i_nj_n}\mapsto
\de_{i_1j_1}\ldots\,\de_{i_nj_n}\cdot(-1)^{\,\bi_1\,+\,\ldots\,+\,\bi_n}\,.
$$
This functional is invariant with respect to the adjoint action
$\operatorname{ad}$
of the Lie superalgebra $\glN$ in $\EndCN^{\,\ot n}$.
Let us now introduce the {\it Capelli element}
$$
C_\la=\tr\ot\id\,(F_\la)\in\UN,
\Tag{7.666666}
$$
see \(7.55555).
This definition of the element $C_\la$ is motivated by the results of
[\hs4\hs,\hs5\hs].

\proclaim{Lemma 4.5}
We have $C_\la\in\ZN$.
\endproclaim

\demo{Proof}
By the equality \(7.2) and by the definition of $F_\la$ we have 
for any $X\in\qN$
$$
\bigl[\,X\,,\,\tr\ot\id\,(F_\la)\,\bigr]=
-\,(\tr\crc\hskip-1pt\operatorname{ad}X)\ot\id\,\bigl(F_\la\bigr)=0
\quad\square
$$
\enddemo

\nt
Note that $\hskip-1pt\ga\,(\ZN)\hskip-1pt\subset\I$.
\hskip-1pt So we get the following corollary to Propositions 
\hskip1pt4.3\hskip1pt,\hskip1pt4.4.

\proclaim{Corollary 4.6}
We have $\ga\,(C_\la)\in\I_\la$.
\endproclaim 

\nt
In the remainder of this section we will
give an explicit formula for the differential
operator $\ga\,(C_\la)\in\PDMN$. 
In particular, we will see that $\ga\,(C_\la)\neq0$ for $\ell_\la\le M,N$. 

Consider the collection \(5.3125) of the basic elements
$h=g\, a_1^{l_1}\ldots a_n^{l_n}$ of the Sergeev algebra $H_n$.
Here $g$ runs through the symmetric group $S_n$ while each of the indices
$l_1\lc l_n$ runs through $0,1$. Let the indices $i_1\lc i_n$
and $b_1\lc b_n$ run through $\pm\,1\lc\pm N$ and $\pm\,1\lc\pm M$
respectively. Put 
$$
I_\la\ =\ 
\sum_g\ 
\sum_{\,l_1\ldots l_n}
\sum_{\,i_1\ldots i_n}
\sum_{\,b_1\ldots b_n}
\chi_\la(h)\cdot
x_{j_nb_n}\ldots x_{j_1b_1}\,\d_{i_1b_1}\ldots\d_{i_nb_n}
\cdot(-1)^e\
\Tag{7.7777777}
$$
where we write $j_s=i_{g(s)}\cdot(-1)^{\,l_s}$ for each $s=1\lc n$ and denote
$$
e\ =\ \sum_{r<s}\ 
(\,\bi_r\,\bb_s+\bj_r\,\bb_s+\bj_r\,\bj_s+\bj_r\,l_s)\ \ +
\!\!\!\!\!\!\!\sum\Sb r<s\\g\1(r)<g\1(s)\endSb\!\!\!\!\!\!\!\!\bi_r\,\bi_s\,
\ +\ \,
\sum_s\ (\bj_s+1)\,l_s\,.
$$

\proclaim{Theorem 4.7}
We have $\ga\,(C_\la)=I_\la$. Here $I_\la\neq0$ if $\ell_\la\le M,N$.
\endproclaim

\demo{Proof}
Consider the $\ZZ$-gradation on the vector space $\PDMN=\PMN\cdot\DMN$
by the degree of the differential operator. By definition any element of
the subspace $\I_\la\subset\PDMN$ is homogeneous of the degree $n$.
Therefore by \(7.666666) due to Corollary 4.6
the element $\ga\,(C_\la)$ coincides with
the leading term of the element
$$
\tr\ot\ga\,\bigl(R_\la\ot1\cdot F_1\ldots F_n\bigr)\in\PDMN.
\Tag{7.7}
$$

Let $Q_\la\in\EndCN^{\,\ot n}$ be the image of the element
$\Chi_\la\in H_n$ under \(3.55555)\,; see \(5.625).
First let us show that
the leading term of \(7.7) coincides with~that~of
$$
\tr\ot\ga\,\bigl(Q_\la\ot1\cdot F_1\ldots F_n\bigr)\in\PDMN.
\Tag{7.8}
$$
Observe that by the definition of the element $F\in\EndCN\ot\UN$ we have
$$
J\ot1\cdot F\cdot J\1\ot1=-\,F.
$$
Therefore for any $s=1\lc n$ we get the equalities in $\PDMN$
$$
\align
\tr&\ot\ga\,\bigl(J_s\,R_\la\,J_s\1\ot1\cdot F_1\ldots F_n\bigr)=
\\
-\,\tr&\ot\ga\,\bigl(J_s\,R_\la\ot1\cdot F_1\ldots F_n\cdot J_s\1\ot1\bigr)=
\tr\ot\ga\,\bigl(R_\la\ot1\cdot F_1\ldots F_n\bigr).
\endalign
$$
Furthermore, for any $s=1\lc n-1$ we have
$$
\align
\tr&\ot\ga\,\bigl(P_{s,s+1}\,R_\la\,P_{s,s+1}\ot1\cdot
F_1\ldots F_s\,F_{s+1}\ldots F_n\bigr)=
\\
\tr&\ot\ga\,\bigl(P_{s,s+1}\,R_\la\ot1\cdot
F_1\ldots F_{s+1}\,F_{s}\ldots F_n\cdot P_{s,s+1}\ot1\bigr)=
\\
\tr&\ot\ga\,\bigl(R_\la\ot1\cdot
F_1\ldots F_{s+1}\,F_{s}\ldots F_n\bigr).
\Tag{7.9}
\endalign
$$
But the elements \(7.7) and \(7.9) of $\PDMN$ have the same leading terms.
So the equlaity of the leading terms in \(7.7) and \(7.8) follows from
Corollary 3.5.

But the leading term of \(7.8) is easy to write down.
Under \(3.55555) for any $g\in S_n$ 
$$
g\1\,\mapsto\sum_{\,i_1\ldots i_n}
E_{i_{g(1)}i_1}\ot\ldots\ot E_{i_{g(n)}i_n}
\cdot\,(-1)^{\,k}
$$
where
$$
k\ \,=\,\ \sum_{r<s}\ \bi_r\,\bigl(\,\bi_s\,+\,\bi_{g(s)}\bigr)\,\,\,+
\!\!\!\!\!\!\!\sum\Sb r<s\\g\1(r)>g\1(s)\endSb\!\!\!\!\!\!\!\!\bi_r\,\bi_s\,.
$$
Further, 
$$
\bigl(a_1^{l_1}\ldots a_n^{l_n}\bigr)\1\,\mapsto
\sum_{\,j_1\ldots j_n}
E_{j_1,\,\ep_1\cdot j_1}\ot\ldots\ot E_{j_n,\,\ep_n\cdot j_n}
\cdot(-1)^{\,l}
$$
where each of the indices $j_1\lc j_n$ runs through $\pm\,1\lc\pm N$\,;
we write $\ep_s=(-1)^{l_s}$ for each $s=1\lc n$ and denote
$$
l\ =\ \sum_{r<s}\ l_r\,l_s\ +\ \sum_s\ l_s\,(\,\bj_s\,+\,1\,).
\nopagebreak
$$
Now the definitions \(5.625) and \(7.1) imply that the leading term
of \(7.8) equals $I_\la$. 

In particular, we have $I_\la\in\I_\la$ by
Corollary 4.6. Suppose that $\ell_\la\le M,N$.
Let us show that $I_\la\neq0$.
Consider the element
$$
I\,=\sum_{\,i_1\ldots i_n}
\sum_{\,b_1\ldots b_n}
x_{i_nb_n}\ldots x_{i_1b_1}\,\d_{i_1b_1}\ldots\d_{i_nb_n}\,\in\,\PDMN
$$
This element is $\UN\ot\UM$-invariant. Due to
[\hs9\hs,\hs Theorem 3\hs] the element 
$I_\la\cdot{\dim U_\la}\,/\,(2^nn!)$
is the projection of $I\in\I$ to the direct summand $\I_\la$ in \(7.0).
This projection cannot be zero since the elements 
$x_{i_1b_1}\ldots x_{i_nb_n}\in\PMN$ span $\PMN_n$
\enddemos

\nt
Consider the {\it canonical} $\ZZ$-filtration on the algebra $\UN$.
It is defined by assigning the degree $1$ to every generator 
$F_{ij}\in\qN$. The corresponding $\ZZ$-graded algebra is
the supersymmetric algebra $\operatorname{S}\,(\qN)$. The subalgebra
$\operatorname{I}\,(\qN)\subset\operatorname{S}\,(\qN)$
of invariants with respect to the adjoint action of $\qN$
corresponds to the centre $\ZN\subset\UN$.

Take
the quotient of the algebra $\operatorname{S}\,(\qN)$ by the ideal
generated by all the elements $F_{ij}$ with $i\neq j$. For each $i=1\lc N$ 
denote by $t_i$ the image in the quotient 
of the element $F_{ii}\in\operatorname{S}\,(\qN)$. Then
the quotient algebra is $\CC\,[\,t_1\lc t_N\,]$. 
Due to [\hs18\hs,\hs Theorem 1\hs] the image in $\CC\,[\,t_1\lc t_N\,]$ 
of the subalgebra
$\operatorname{I}\,(\qN)\subset\operatorname{S}\,(\qN)$
is generated by all the power sums of odd degree
$\,t_1^{\,2r+1}+\ldots+t_N^{\,2r+1}\,$, $r\ge0$.
We will desribe the
image $T_\la\in\CC\,[\,t_1\lc t_N\,]$ of the element
from $\operatorname{I}\,(\qN)\subset\operatorname{S}\,(\qN)$
corresponding to the central element $C_\la\in\UN$.
In particular, we will see that $C_\la\neq0$.
For description of the eigenvalue of $C_\la\in\UN$
in the irreducible module $V_\mu$ where the diagram
$\mu$ does contain $\la$,~\text{see [\hs7\hs].}

A {\it shifted Young tableau} of shape $\la$ is any
bijective filling of the boxes of diagram $\la$ with the numbers
$1\lc n$ such that in every row and column the numbers increase from
the left to the right and from the top to the bottom respectively.
Let $n_\la$ stand for the total number of shifted Young tableaux of
the shape $\la$.
Let $Q_\la(\,t_1\lc t_N\,;-\,1\,)$ be the Schur 
$Q$-polynomial symmetric in $t_1\lc t_N$; see [\hs{21}\hs].

\proclaim{Proposition 4.8}
We have $T_\la(\,t_1\lc t_N\,)=Q_\la(\,t_1\lc t_N\,;-\,1\,)\cdot n!/n_\la\,$.
\endproclaim

\demo{Proof}
Regard $t_1\lc t_N$ as complex variables. Put $t_{-i}=t_i$
for each $i=1\lc N$. Denote 
$$
T\,=\,\sum_i\ t_i\,E_{ii}\cdot(-1)^\bi\in\EndCN
$$
where the index $i$ runs through $\pm1\lc\pm N$. Then by the definition
\(7.666666) we have
$T_\la=\tr\,\bigl(\hskip1pt Q_\la\cdot T^{\,\ot n}\hskip1pt\bigr)$;
see also the proof of Theorem 4.7. Now Proposition 4.8 can be obtained
from [\hs9\hs,\hs Section 2.2\hs] and [\hs21\hs,\hs Example III.7.8\hs]
\enddemos

\nt
Due to Proposition 4.8 the elements $C_\la$ where the diagram $\la$
has only one row, generate the centre
of enveloping algebra $\UN$. So the equality $\ga\,(C_\la)=I_\la$ for
$\la=(n)$ may be regarded as an analogue for $\qN$ of the classical
\text{Capelli identity [\hs1\hs].} Note that for $\la=(n)$
the element $\Psi_\la\in\Sen$ is determined by the 
formula \(1.0101).


\kern15pt
\centerline{\bf Acknowledgements}\section{\,}
\kern-20pt

\nt
The author has been supported by the EPSRC Advanced Research Fellowship.
He was also supported by the Erwin Schr\"odinger International Institute
in Vienna.


\newpage
\vbox{
\centerline{\bf References}\section{\,}
\kern-22pt

\itemitem{[1]}
{A. Capelli},
{\it Sur les op\'erations dans la th\'eorie des formes alg\'ebriques},
{Math. Ann.}
{\bf 37}
(1890),
1--37.

\itemitem{[2]}
{H. Weyl},
{\it Classical groups, their invariants and representations},
Princeton University Press, Princeton, 1946.

\itemitem{[3]}
{S. Sahi},
{\it The spectrum of certain invariant differential operators associated 
to a Hermitian symmetric space,}
in
{\lq\lq\ts Lie Theory and Geometry\ts\rq\rq}, 
{Progress in Math.}
{\bf 123},
Birkh\"auser,
Boston,
1994,
pp. 569--576.

\itemitem{[4]}
{A. Okounkov},
{\it Quantum immanants and higher Capelli identities},
{Transformation Groups}
{\bf 1}
(1996),
99--126.

\itemitem{[5]}
{M. Nazarov},
{\it Young's symmetrizers for projective representations of
the symmetric group}, to appear in Adv. Math.

\itemitem{[6]}
V. Kac,
{\it Lie superalgebras},
Adv. Math.
{\bf 26}
(1977),
8--96.

\itemitem{[7]}
V. Ivanov, 
{\it Dimension of a skew shifted Young diagram and projective
representations of the infinite symmetric group},
to appear in J. Math. Sciences.

\itemitem{[8]}
A. Borodin and N. Rozhkovskaya,
{\it On a superanalogue of the Schur--Weyl duality},
Preprint ESI 
{\bf 246}
(1995),
1-15.

\itemitem{[9]}
{A. Sergeev,}
{\it The tensor algebra of the identity representation 
\text{as a module over}
the Lie superalgebras $GL(n,m)$ and $Q(n)$},
{Math.\,Sbornik}
{\bf 51}
(1985),
\text{419--427.}

\itemitem{[10]}
{I. Cherednik},
{\it On special bases of irreducible finite-dimensional representations
of the degenerate affine Hecke algebra},
{Funct. Anal. Appl.}
{\bf 20}
(1986),
87--89.

\itemitem{[11]}
{M. Nazarov},
{\it Young's symmetrizers for projective representations of
the symmetric group}, to appear in Adv. Math.

\itemitem{[12]}
{A. Okounkov},
{\it Young basis, Wick formula, and higher Capelli identities},
to appear in {Int. Math. Research Notes};
{q-alg/9602027}.

\itemitem{[13]}
{A. Molev},
{\it Factorial supersymmetric Schur functions and super Capelli identities},
to appear in Amer. Math. Soc. Transl.;
{q-alg/9606008}.

\itemitem{[14]} 
{M. Nazarov},
{\it Yangian of the queer Lie superalgebra},
in \lq\lq\,Symposium in Representation Theory\,\rq\rq,
Yamagata, 1992,
pp. 8--15.

\itemitem{[15]}
M. Jimbo, A. Kuniba, T. Miwa and M. Okado,
{\it The $A_n^{(1)}$ face models},
Comm. Math. Phys.
{\bf 119}
(1988),
543--565.

\itemitem{[16]}
{A. Young,}
{\it On quantitative substitutional analysis I} and {\it II},
{Proc. London Math. Soc.}
{\bf 33}
(1901),
97--146
{and}
{\bf 34}
(1902),
361--397.

\itemitem{[17]}
{A. Jucys,}
{\it Symmetric polynomials and the center of the symmetric group ring\,},
{Rep.\, Math.\, Phys.}
{\bf 5}
(1974),
107--112.

\itemitem{[18]}
{A. Sergeev,}
{\it The centre of enveloping algebra of the the Lie superalgebra $Q(n)$},
{Lett. Math. Phys.}
{\bf 7}
(1983),
177--179.

\itemitem{[19]}
{R. Howe},
{\it Remarks on classical invariant theory},
{Trans. Amer. Math. Soc.}
{\bf 313}
(1989),
539--570.

\itemitem{[20]}
{M. Duflo},
{\it Sur la classification des id\'eaux primitifs dans
l'alg\graveaccent{e}bre
enveloppante d'une alg\graveaccent{e}bre de Lie semi-simple},
{Ann. of Math.}
{\bf 105}
(1977),
107--120.

\itemitem{[21]}
{I. Macdonald},
{\it Symmetric functions abnd Hall polynomials},
Clarendon Press, Oxford, 1979.


\kern-2.53339pt
\centerline{\hbox to 3.2cm{\hrulefill}}
\kern9pt
\centerline{Mathematics Department, University of York, York YO1 5DD, England}
\centerline{E-mail: {mln1@york.ac.uk}}
}


\bye